# FOCK-SOBOLEV SPACES OF FRACTIONAL ORDER

HONG RAE CHO, BOO RIM CHOE, AND HYUNGWOON KOO

ABSTRACT. For the full range of index $0 < p \leq \infty$, real weight $\alpha$ and real Sobolev order $s$, two types of weighted Fock-Sobolev spaces over $\mathbb{C}^n$, $F_{\alpha,s}^p$ and $\widetilde{F}_{\alpha,s}^p$, are introduced through fractional differentiation and through fractional integration, respectively. We show that they are the same with equivalent norms and, furthermore, that they are identified with the weighted Fock space $F_{\alpha-sp,0}^p$ for the full range of parameters. So, the study on the weighted Fock-Sobolev spaces is reduced to that on the weighted Fock spaces. We describe explicitly the reproducing kernels for the weighted Fock spaces and then establish the boundedness of integral operators induced by the reproducing kernels. We also identify dual spaces, obtain complex interpolation result and characterize Carleson measures.

## 1. INTRODUCTION

Function theoretic and also operator theoretic properties of Fock space have been studied widely for the last several years. We refer the reader to [8] and [11] for more recent and systematic treatment of Fock spaces. Recently Cho and Zhu [5] studied Fock-Sobolev spaces of positive integer order over the multi-dimensional complex spaces. The purpose of the current paper is to extend the notion of the Fock-Sobolev spaces to the case of fractional orders allowed to be any real number. Most of our results, even when restricted to the case of positive integer orders, contain the results in [5] as special cases.

Throughout the paper $n$ is a fixed positive integer, reserved for the dimension of the underlying multi-dimensional complex space. We write $dV$ for the volume measure on the complex $n$-space $\mathbb{C}^n$ normalized so that $\int_{\mathbb{C}^n} e^{-|z|^2} \, dV(z) = 1$. Also, we write $z \cdot \overline{w}$ for the Hermitian inner product of $z, w \in \mathbb{C}^n$ and let $|z| = (z \cdot \overline{z})^{1/2}$. More explicitly,

$$z \cdot \overline{w} = \sum_{j=1}^n z_j \overline{w}_j, \qquad |z| = \left( \sum_{j=1}^n |z_j|^2 \right)^{1/2}$$

Date: July 9, 2013.

2010 *Mathematics Subject Classification.* Primary 32A37; Secondary 30H20.

H. Cho was supported by the National Research Foundation of Korea(NRF) grant funded by the Korea government(MEST) (NRF-2011-0013740) and B. Choe was supported by Basic Science Research Program through the National Research Foundation of Korea(NRF) funded by the Ministry of Education, Science and Technology(2013R1A1A2004736).

where $z_j$ denotes the $j$-th component of a typical point $z \in \mathbb{C}^n$ so that $z = (z_1, \ldots, z_n)$.

It will turn out that polynomially growing/decaying weights quite naturally come into play in the study of our Fock-Sobolev spaces of fractional order. So, we first introduce such weighted Fock spaces. Given $\alpha$ real we put

$$dV_\alpha(z) = \frac{dV(z)}{(1 + |z|)^\alpha}. \tag{1.1}$$

Now, for $0 < p < \infty$, we denote by $L_\alpha^p = L_\alpha^p(\mathbb{C}^n)$ the space of Lebesgue measurable functions $\psi$ on $\mathbb{C}^n$ such that the norm

$$\|\psi\|_{L_\alpha^p} := \left\{ \int_{\mathbb{C}^n} \left| \psi(z) e^{-\frac{1}{2}|z|^2} \right|^p \, dV_\alpha(z) \right\}^{1/p}$$

is finite; here, we are abusing the term "norm" for $0 < p < 1$ only for convenience. For $p = \infty$, we denote by $L_\alpha^\infty = L_\alpha^\infty(\mathbb{C}^n)$ the space of Lebesgue measurable functions $\psi$ on $\mathbb{C}^n$ such that the norm

$$\|\psi\|_{L_\alpha^\infty} := \operatorname{esssup} \left\{ \frac{|\psi(z)| e^{-\frac{1}{2}|z|^2}}{(1 + |z|)^\alpha} : z \in \mathbb{C}^n \right\} \tag{1.2}$$

is finite.

Now, for $\alpha$ real and $0 < p \le \infty$, we define

$$F_\alpha^p := L_\alpha^p \cap H(\mathbb{C}^n)$$

where $H(\mathbb{C}^n)$ denotes the class of entire functions on $\mathbb{C}^n$. Of course, we regard $F_\alpha^p$ as a subspace of $L_\alpha^p$. The space $F_\alpha^p$ is closed in $L_\alpha^p$ and thus is a Banach space when $1 \le p \le \infty$. In particular, $F_\alpha^2$ is a Hilbert space. Also, for $0 < p < 1$, the space $F_\alpha^p$ is a complete metric space under the translation-invariant metric $(f, g) \mapsto \|f - g\|_{F_\alpha^p}^p$; see the remark at the end of Section 2.

We write $\|f\|_{F_\alpha^p} := \|f\|_{L_\alpha^p}$ for $f \in H(\mathbb{C}^n)$ in order to emphasize that $f$ is holomorphic. Also, we write $F^p = F_0^p$ when $\alpha = 0$. The space $F^p$ is often called under the various different names such as Fock space, Bargmann space, Segal-Bargmann space, and so on. We call it Fock space for no particular reason. Naturally we call the space $F_\alpha^p$ a weighted Fock space.

We now introduce two different types of weighted Fock-Sobolev spaces of fractional order: one in terms of fractional differentiation operator $\mathcal{R}^s$ and the other in terms of fractional integration operator $\widetilde{\mathcal{R}}^s$. The precise definitions of $\mathcal{R}^s$ and $\widetilde{\mathcal{R}}^s$ are given in Section 3. Given any real number $\alpha$ and $s$, the first type of weighted Fock-Sobolev space $F_{\alpha,s}^p$ is defined to be the space of all $f \in H(\mathbb{C}^n)$ such that $\mathcal{R}^s f \in L_\alpha^p$. The second type of weighted Fock-Sobolev space $\widetilde{F}_{\alpha,s}^p$ is defined similarly with $\widetilde{\mathcal{R}}^s$ in place of $\mathcal{R}^s$. The precise norms on these weighted Fock-Sobolev spaces are given in Section 3. We refer to [4], [6] and [7] for other Sobolev spaces of similar type.

Our result (Theorem 4.2) shows that two notions of weighted Fock-Sobolev spaces coincide and that they can be realized as a weighted Fock space: for any $\alpha$ and $s$ real,

$$F_{\alpha,s}^p = F_{\alpha-sp}^p = \widetilde{F}_{\alpha,s}^p \qquad \text{for} \quad 0 < p < \infty \tag{1.3}$$

and

$$F_{\alpha,s}^\infty = F_{\alpha-s}^\infty = \widetilde{F}_{\alpha,s}^\infty \qquad \text{for} \quad p = \infty \tag{1.4}$$

with equivalent norms. Section 4 is devoted to the proof of these characterizations. Note that the most natural definition of the weighted Fock-Sobolev space of positive integer order might be the one in terms of full derivatives. That turns out to be actually the case as a consequence of the first equalities in (1.3) and (1.4); see Corollary 4.4. For the unweighted case such a characterization in terms of full derivatives has been already noticed in [4] for $p = 2$ and [5] for general $0 < p < \infty$. Also, the result (1.3) is quite reminiscent of what have been known for the weighted Bergman-Sobolev spaces $A_{\alpha,s}^p(\mathbb{B}_n)$ over the unit ball $\mathbb{B}_n$ of $\mathbb{C}^n$:

$$A_{\alpha,s}^p(\mathbb{B}_n) = A_{\alpha-sp,0}^p(\mathbb{B}_n) = A_{0,s-\alpha/p}^p(\mathbb{B}_n)$$

with equivalent norms. In this ball case, however, the weight $(1-|z|^2)^\alpha$ is restricted to $\alpha > -1$, the order $s$ of fractional differentiation is restricted to $s \geq 0$ and the index $p$ is restricted to $\alpha - sp > -1$; see [2] and [9].

As key preliminary steps towards (1.3) and (1.4), we describe how the fractional differentiation/integration act on the weighted Fock spaces (Theorem 3.13). In the course of the proof we obtain integral representations for fractional differentiation/integration and use them to establish pointwise size estimates of the fractional derivative/integral of the well-known Fock kernel $e^{z \cdot \overline{w}}$. These results are proved in Section 3.

Having characterizations (1.3) and (1.4), we may focus on weighted Fock spaces in order to study properties of weighted Fock-Sobolev spaces. As is easily seen in Section 4, the weighted Fock space $F_\alpha^2$ is a reproducing kernel Hilbert space. For example, the aforementioned Fock kernel is the reproducing kernel for the unweighted Fock space $F^2$. We obtain an explicit descriptions (Theorem 4.5) of the reproducing kernels.

As applications we derive some fundamental properties of the weighted Fock-Sobolev spaces such as:

- Reproducing operator;
- Dual space;
- Complex interpolation;
- Carleson measure.

These results are proved in Section 5.

*Constants.* In this paper we use the same letter $C$ to denote various positive constants which may vary at each occurrence but do not depend on the essential parameters. Variables indicating the dependency of constants $C$ will be often specified in parenthesis. For nonnegative quantities $X$ and $Y$ the notation $X \lesssim Y$ or $Y \gtrsim X$ means $X \leq CY$ for some inessential constant $C$. Similarly, we write $X \approx Y$ if both $X \lesssim Y$ and $Y \lesssim X$ hold.

## 2. Some basic properties

basic

In this section we observe two basic properties for the weighted Fock spaces. One is the growth estimate of weighted Fock functions and the other is the density of holomorphic polynomials.

mvplem

**Lemma 2.1.** *Given $a, t > 0$ and $\alpha$ real, there is a constant $C = C(a, t, \alpha) > 0$ such that*

$$\frac{|f(z)|^p}{e^{a|z|^2}(1+|z|)^\alpha} \leq C \int_{|w-z|<t} |f(w)|^p e^{-a|w|^2} \, dV_\alpha(w), \qquad z \in \mathbb{C}^n,$$

*for $0 < p < \infty$ and $f \in H(\mathbb{C}^n)$.*

*Proof.* We first mention an elementary inequality

$$\left(\frac{1+|z|}{1+|w|}\right)^\alpha \leq (1+|z-w|)^{|\alpha|} \tag{2.1}$$

element

valid for any $\alpha$ real and $z, w \in \mathbb{C}^n$. To see this, note $1 + |z| \leq 1 + |w| + |z-w|$ and therefore

$$\frac{1+|z|}{1+|w|} \leq 1 + |z-w|$$

for any $z, w \in \mathbb{C}^n$.

Let $a, t > 0$ and $\alpha$ be a real number. Let $0 < p < \infty$ and $f \in H(\mathbb{C}^n)$. Fix $z \in \mathbb{C}^n$. We have by subharmonicity of the function $w \mapsto |f(z+w)e^{-2aw\cdot\overline{z}/p}|^p$

$$|f(z)|^p \lesssim \frac{1}{t^{2n}} \int_{|w|<t} \left|f(z+w)e^{-2aw\cdot\overline{z}/p}\right|^p \, dV(w)$$

$$\leq \frac{e^{at^2}}{t^{2n}} \int_{|w|<t} \left|f(z+w)e^{-2aw\cdot\overline{z}/p}\right|^p e^{-a|w|^2} \, dV(w)$$

$$= \frac{e^{at^2}}{t^{2n}} e^{a|z|^2} \int_{|w-z|<t} |f(w)|^p e^{-a|w|^2} \, dV(w).$$

Note by (2.1)

$$1 < (1+t)^{|\alpha|} \left(\frac{1+|z|}{1+|w|}\right)^\alpha$$

for $|w-z| < t$. Combining these observations, we conclude the asserted inequality. $\square$

In what follows we use the standard multi-index notation. Namely, given an $n$-tuple $\gamma = (\gamma_1, \ldots, \gamma_n)$ of nonnegative integers, $|\gamma| = \sum_{j=1}^n \gamma_j$ and $\partial^\gamma = \partial_1^{\gamma_1} \cdots \partial_n^{\gamma_n}$, etc., where $\partial_j = \partial/\partial z_j$.

cauchy

**Proposition 2.2.** *Given $0 < p \leq \infty$, $\alpha$ real and a multi-index $\gamma$, there is a constant $C = C(p, \alpha, \gamma) > 0$ such that*

$$|\partial^\gamma f(z)| \leq C e^{\frac{|z|^2}{2}} (1+|z|)^{\frac{\alpha}{p}+|\gamma|} \|f\|_{F_\alpha^p}, \qquad 0 < p < \infty$$

*and*

$$|\partial^\gamma f(z)| \le C e^{\frac{|z|^2}{2}} (1 + |z|)^{\alpha + |\gamma|} \|f\|_{F_\alpha^\infty}$$

*for $z \in \mathbb{C}^n$ and $f \in H(\mathbb{C}^n)$.*

*Proof.* Fix $\alpha$ real and consider the case $0 < p < \infty$. The case $\gamma = 0$ is an immediate consequence of Lemma 2.1 (with $a = p/2$). Let $f \in H(\mathbb{C}^n)$ and $z \in \mathbf{C}^n$. We may assume $|z| \ge 1$. Given a multi-index $\gamma$, applying the Cauchy estimates on the ball with center $z$ and radius $1/|z|$, we have by the maximum modulus theorem and Lemma 2.1

$$\begin{aligned}
|\partial^\gamma f(z)| &\lesssim |z|^{|\gamma|} \max_{|w-z|=1/|z|} |f(w)| \\
&\le |z|^{|\gamma|} \max_{|w|=|z|+1/|z|} |f(w)| \\
&\lesssim e^{(|z|+1/|z|)^2/2} |z|^{|\gamma|} (1 + |z| + 1/|z|)^{\frac{\alpha}{p}} \|f\|_{F_\alpha^p}.
\end{aligned}$$

Meanwhile, since $|z| \ge 1$, we have

$$e^{(|z|+1/|z|)^2/2} = e^{(|z|^2 + 2 + 1/|z|^2)/2} \approx e^{|z|^2/2}$$

and

$$|z|^{|\gamma|} (1 + |z| + 1/|z|)^{\frac{\alpha}{p}} \approx (1 + |z|)^{\frac{\alpha}{p} + |\gamma|}.$$

Thus we conclude the asserted estimate for $p$ finite.

When $p = \infty$, note that the case $\gamma = 0$ holds by definition of $F_\alpha^\infty$. So, we have the asserted estimate by the same argument. The proof is complete. $\square$

As one may quite naturally expect, holomorphic polynomials form a dense subset in any weighted Fock space with $p$ finite. To see it we first note a basic fact:

$$\lim_{r \to 1^-} \|f_r - f\|_{F_\alpha^p} = 0 \qquad (2.2) \quad \boxed{\texttt{dilation}}$$

where $f_r(z) = f(rz)$ for $0 < r < 1$. This follows from the fact $\|f_r\|_{F_\alpha^p}^p \to \|f\|_{F_\alpha^p}^p$ as $r \to 1^-$, which can be easily verified via an elementary change-of-variable and the dominated convergence theorem.

Note that (2.2) does not extend to the case $p = \infty$. In conjunction with this observation, we introduce a subspace of $F_\alpha^\infty$ that enjoys the property (2.2). Given $\alpha$ real, let $F_\alpha^{\infty,0}$ be the space consisting of all $f \in F_\alpha^\infty$ such that

$$\lim_{|z| \to \infty} \frac{|f(z)| e^{-\frac{|z|^2}{2}}}{(1 + |z|)^\alpha} = 0. \qquad (2.3) \quad \boxed{\texttt{littledef}}$$

It is easily checked that $F_\alpha^{\infty,0}$ is a closed subspace of $F_\alpha^\infty$. Also, for $f \in F_\alpha^\infty$, we have $f \in F_\alpha^{\infty,0}$ if and only if (2.2) with $p = \infty$ holds.

$\boxed{\texttt{dense}}$ **Proposition 2.3.** *Given $\alpha$ real, the set of all holomorphic polynomials is dense in $F_\alpha^{\infty,0}$ and $F_\alpha^p$ for any $0 < p < \infty$.*

*Proof.* We modify the proof of [11, Proposition 2.9] where the one-variable version of the unweighted case is treated. Fix $\alpha$ real.

We first consider the case $0 < p < \infty$. Let $f \in F_\alpha^p$. By (2.2) it suffices to show that the homogeneous expansion of $f_r$ converges in $F_\alpha^p$ for each $0 < r < 1$. Namely, using the homogeneous expansion $f = \sum_{k=0}^\infty f_k$ where $f_k$ is a homogeneous polynomial of degree $k$, it is enough to prove

$$\left\| \sum_{k=N}^\infty r^k f_k \right\|_{F_\alpha^p} \to 0 \qquad (r : \text{fixed}) \tag{2.4}$$

as $N \to \infty$.

In order to establish (2.4) we need to estimate the size of Taylor coefficients and the norms of monomials. To estimate the size of Taylor coefficients, we note for a given multi-index $\nu$ by the Cauchy integral formula over the unit polydisk

$$t^\nu \partial^\nu f(0) = \frac{\nu!}{(2\pi i)^n} \int_{|\zeta_1|=1} \cdots \int_{|\zeta_n|=1} \frac{f(t_1\zeta_1, \ldots, t_n\zeta_n)}{\zeta_j^{\nu_j+1}} \, d\zeta_1 \cdots d\zeta_n$$

for any $t = (t_1, \ldots, t_n)$ where $t_j > 0$ if $\nu_j > 0$ and $t_j = 0$ otherwise. Since $|(t_1\zeta_1, \ldots, t_n\zeta_n)| = |t|$, the above and Proposition 2.2 yield

$$\left| \frac{\partial^\nu f(0)}{\nu!} \right| \lesssim \frac{1}{t^\nu} e^{\frac{|t|^2}{2}} (1 + |t|)^{\frac{\alpha}{p}} \|f\|_{F_\alpha^p}.$$

So, choosing $t_j = \sqrt{\nu_j}$ when $\nu_j > 0$, we have

$$\left| \frac{\partial^\nu f(0)}{\nu!} \right| \lesssim e^{\frac{|\nu|}{2}} |\nu|^{\frac{\alpha}{2p}} \left( \prod_{j=1}^n \nu_j^{-\frac{\nu_j}{2}} \right) \|f\|_{F_\alpha^p} \tag{2.5}$$

for $|\nu|$ large where $\nu_j^{-\frac{\nu_j}{2}}$ is understood to be 1 when $\nu_j = 0$.

To estimate the norms of monomials, we note

$$\int_S |\zeta^\nu|^p \, d\sigma(\zeta) = \frac{\Gamma(n)}{\Gamma\left(\frac{p}{2}|\nu| + n\right)} \prod_{j=1}^n \Gamma\left(\frac{p}{2}\nu_j + 1\right);$$

to see this one may easily modify the proof of [9, Lemma 1.11] where the case $p = 2$ is proved. Thus, integrating in polar coordinates, we obtain

$$\|z^\nu\|_{F_\alpha^p}^p = c_n \frac{\prod_{j=1}^n \Gamma(\frac{p}{2}\nu_j + 1)}{\Gamma(\frac{p}{2}|\nu| + n)} \int_0^\infty \frac{e^{-\frac{p}{2}t^2} t^{p|\nu|+2n-1}}{(1+t)^\alpha} \, dt$$

for some dimensional constant $c_n$. Meanwhile, we have for $|\nu| > (2n - \alpha)/p$

$$\int_0^\infty \frac{e^{-\frac{p}{2}t^2} t^{p|\nu|+2n-1}}{(1+t)^\alpha} \, dt \lesssim \int_0^\infty e^{-\frac{p}{2}t^2} t^{p|\nu|+2n-1-\alpha} \, dt$$
$$= \frac{1}{p} \left( \frac{2}{p} \right)^{\frac{p}{2}|\nu| - \frac{\alpha}{2} + n - 1} \Gamma\left( \frac{p}{2}|\nu| - \frac{\alpha}{2} + n \right).$$

So far, we have

$$\|z^\nu\|_{F_\alpha^p}^p \lesssim \left(\frac{2}{p}\right)^{p|\nu|/2} \frac{\Gamma\left(\frac{p}{2}|\nu| - \frac{\alpha}{2} + n\right)}{\Gamma(\frac{p}{2}|\nu| + n)} \prod_{j=1}^n \Gamma\left(\frac{p}{2}\nu_j + 1\right) \qquad (2.6) \quad \boxed{\texttt{sofar}}$$

for $|\nu|$ large. Since

$$\frac{\Gamma(\frac{p}{2}|\nu| - \frac{\alpha}{2} + n)}{\Gamma(\frac{p}{2}|\nu| + n)} \approx \left(\frac{p}{2}|\nu| - \frac{\alpha}{2} + n\right)^{-\frac{\alpha}{2}} \approx |\nu|^{-\frac{\alpha}{2}}$$

by Stirling's formula, we obtain from (2.6)

$$\|z^\nu\|_{F_\alpha^p} \lesssim \left(\frac{2}{p}\right)^{\frac{|\nu|}{2}} |\nu|^{-\frac{\alpha}{2p}} \prod_{j=1}^n \Gamma^{\frac{1}{p}}\left(\frac{p}{2}\nu_j + 1\right) \qquad (2.7) \quad \boxed{\texttt{mononorm}}$$

for $|\nu|$ large. Meanwhile, we have by Stirling's formula

$$\prod_{j=1}^n \Gamma\left(\frac{p}{2}\nu_j + 1\right) \lesssim \prod_{j=1}^n \left[\nu_j^{\frac{p}{2}\nu_j + \frac{1}{2}} \left(\frac{p}{2}\right)^{\frac{p}{2}\nu_j + \frac{1}{2}} e^{-\frac{p}{2}\nu_j}\right]$$

$$= \left(\frac{p}{2}\right)^{\frac{p}{2}|\nu| + \frac{n}{2}} e^{-\frac{p}{2}|\nu|} \prod_{j=1}^n \nu_j^{\frac{p}{2}\nu_j + \frac{1}{2}}$$

so that

$$\prod_{j=1}^n \Gamma^{\frac{1}{p}}\left(\frac{p}{2}\nu_j + 1\right) \lesssim \left(\frac{p}{2}\right)^{\frac{|\nu|}{2}} e^{-\frac{|\nu|}{2}} |\nu|^{\frac{n}{2p}} \prod_{j=1}^n \nu_j^{\frac{\nu_j}{2}}$$

for $|\nu|$ large. Thus we have

$$\|z^\nu\|_{F_\alpha^p} \lesssim |\nu|^{\frac{1}{2p}(n-\alpha)} e^{-\frac{|\nu|}{2}} \prod_{j=1}^n \nu_j^{\frac{\nu_j}{2}} \qquad (2.8) \quad \boxed{\texttt{znorm}}$$

for $|\nu|$ large.

Consequently, we have by (2.5) and (2.8)

$$\left|\frac{\partial^\nu f(0)}{\nu!}\right| \|z^\nu\|_{F_\alpha^p} \lesssim |\nu|^{\frac{n}{2p}} \|f\|_{F_\alpha^p} \qquad (2.9) \quad \boxed{\texttt{nuterm}}$$

for $|\nu|$ large and thus

$$\|f_k\|_{F_\alpha^p} \lesssim \sum_{|\nu|=k} \left|\frac{\partial^\nu f(0)}{\nu!}\right| \|z^\nu\|_{F_\alpha^p} \lesssim k^{\frac{n}{2p}}(1+k)^n \|f\|_{F_\alpha^p} \approx k^{\frac{n}{2p}+n} \|f\|_{F_\alpha^p}$$

for $k$ large. Now, for $1 \le p < \infty$, we have

$$\left\|\sum_{k=N}^\infty r^k f_k\right\|_{F_\alpha^p} \le \sum_{k=N}^\infty r^k \|f_k\|_{F_\alpha^p} \lesssim \|f\|_{F_\alpha^p} \sum_{k=N}^\infty r^k k^{\frac{n}{2p}+n} \to 0$$

as $N \to \infty$. On the other hand, for $0 < p < 1$, we have

$$\left\| \sum_{k=N}^{\infty} r^k f_k \right\|_{F_\alpha^p}^p \leq \sum_{k=N}^{\infty} r^{pk} \|f_k\|_{F_\alpha^p}^p \lesssim \|f\|_{F_\alpha^p}^p \sum_{k=N}^{\infty} r^{pk} k^{\frac{n}{2}+np} \to 0$$

as $N \to \infty$. This completes the proof of (2.4) and thus the proof for the case $0 < p < \infty$.

Now, we consider the case $p = \infty$. We claim that there is a constant $C = C(\alpha) > 0$ such that

$$\left| \frac{\partial^\nu f(0)}{\nu!} \right| \|z^\nu\|_{F_\alpha^\infty} \leq C \|f\|_{F_\alpha^\infty} \qquad (2.10) \quad \boxed{\texttt{ffa}}$$

for all multi-indices $\nu$ and $f \in F_\alpha^\infty$. With this granted, we see that (2.4) with $p = \infty$ remains valid for $f \in F_\alpha^\infty$ and hence deduce from (2.2) (with $p = \infty$ valid for functions in $F_\alpha^{\infty,0}$) that holomorphic polynomials forms a dense subset in $F_\alpha^{\infty,0}$.

It remains to show (2.10). Let $f \in F_\alpha^\infty$. Note by a trivial modification of the proof of (2.5)

$$\left| \frac{\partial^\nu f(0)}{\nu!} \right| \lesssim e^{\frac{|\nu|}{2}} |\nu|^{\frac{\alpha}{2}} \left( \prod_{j=1}^{n} \nu_j^{-\frac{\nu_j}{2}} \right) \|f\|_{F_\alpha^\infty}$$

for $|\nu|$ large. On the other hand, since

$$\|z^\nu\|_{F_\alpha^\infty} \approx \sup_{|z| \geq 1} |z^\nu| |z|^{-\alpha} e^{-\frac{1}{2}|z|^2} = \left( \sup_{|\zeta|=1} |\zeta^\nu| \right) \left( \sup_{t \geq 1} t^{|\nu|-\alpha} e^{-\frac{t^2}{2}} \right),$$

an elementary calculation yields

$$\|z^\nu\|_{F_\alpha^\infty} \approx |\nu|^{-|\nu|/2} \left( \prod_{j=1}^{n} \nu_j^{\frac{\nu_j}{2}} \right) (|\nu| - \alpha)^{\frac{|\nu|-\alpha}{2}} e^{-\frac{|\nu|-\alpha}{2}}$$

for $|\nu|$ large. It follows that

$$\left| \frac{\partial^\nu f(0)}{\nu!} \right| \|z^\nu\|_{F_\alpha^\infty} \lesssim \left( 1 - \frac{\alpha}{|\nu|} \right)^{\frac{|\nu|-\alpha}{2}} e^{\frac{\alpha}{2}} \|f\|_{F_\alpha^\infty} \to \|f\|_{F_\alpha^\infty}$$

as $|\nu| \to \infty$. So, (2.10) holds, as required. The proof is complete. $\qquad \square$

We now close the section with the following remark for $0 < p \leq \infty$ and $\alpha$ real.

*Remark.* (1) As a consequence of Proposition 2.2 we see that the convergence in the weighted Fock spaces implies the uniform convergence on compact sets. Accordingly, the space $F_\alpha^p$ is closed in $L_\alpha^p$.

(2) When $p < \infty$, in addition to Proposition 2.2, we also have

$$\lim_{|z| \to \infty} \frac{|\partial^\gamma f(z)| e^{-\frac{|z|^2}{2}}}{(1+|z|)^{\frac{\alpha}{p}+|\gamma|}} = 0 \qquad (2.11) \quad \boxed{\texttt{little}}$$

for any multi-index $\gamma$ and $f \in F_\alpha^p$. This can be easily verified by Proposition 2.2 and (2.2).

(3) We mention an estimate to be used later. Given a nonnegative integer $m$, there is a constant $C = C(p, \alpha, m) > 0$ such that

$$\sup_{|w| \leq 1} \frac{|f(w) - f_m(w)|}{|w|^{m+1}} \leq C \|f\|_{F_\alpha^p} \qquad (2.12) \quad \boxed{\texttt{falpha}}$$

for $f \in F_\alpha^p$ where $f_m$ is the Taylor polynomial of $f$ degree $m$. To see this one may apply Proposition 2.2 together with Taylor's formula.

$\boxed{\texttt{fractional}}$

## 3. FRACTIONAL DIFFERENTIATION/INTEGRATION

In this section we define the fractional differentiation/integration and then show how they act on the weighted Fock spaces.

Given $s$ real and $f \in H(\mathbb{C}^n)$ with homogeneous expansion

$$f = \sum_{k=0}^{\infty} f_k \qquad (3.1) \quad \boxed{\texttt{homoexp}}$$

where $f_k$ is a homogeneous polynomial of degree $k$, we define the fractional derivative $\mathcal{D}^s f$ of order $s$ as follows:

$$\mathcal{D}^s f = \begin{cases} \displaystyle\sum_{k=0}^{\infty} \frac{\Gamma(n+s+k)}{\Gamma(n+k)} f_k & \text{if} \quad s \geq 0 \\[2ex] \displaystyle\sum_{k>|s|} \frac{\Gamma(n+s+k)}{\Gamma(n+k)} f_k & \text{if} \quad s < 0. \end{cases} \qquad (3.2) \quad \boxed{\texttt{fracdif}}$$

We remark that our definition of $\mathcal{D}^s f$ is slightly different from the usual ones on the unit ball which is defined as $\sum k^s f_k$ or $\sum (1+k)^s f_k$, but they are asymptotically the same in the sense that $\frac{\Gamma(n+s+k)}{\Gamma(n+k)} \sim k^s$ as $k \to \infty$ by Stirling's formula.

Next, we define the fractional integral $\mathcal{I}^s f$ of order $s$ as follows:

$$\mathcal{I}^s f = \begin{cases} \displaystyle\sum_{k=0}^{\infty} \frac{\Gamma(n+k)}{\Gamma(n+s+k)} f_k & \text{if} \quad s \geq 0 \\[2ex] \displaystyle\sum_{k>|s|} \frac{\Gamma(n+k)}{\Gamma(n+s+k)} f_k & \text{if} \quad s < 0. \end{cases} \qquad (3.3) \quad \boxed{\texttt{fracint}}$$

It is elementary to check that the series above converge uniformly on compact sets and thus $\mathcal{D}^s f$ and $\mathcal{I}^s f$ are again entire functions. Note that $\mathcal{D}^s$ is essentially the inverse operator of $\mathcal{I}^s$, and vice versa.

We first establish pointwise size estimates for the fractional derivatives/integrals of the Fock kernel given by

$$K_w(z) = K(z, w) := e^{z \cdot \overline{w}}$$

for $z, w \in \mathbb{C}^n$. As is well known, this Fock kernel has the reproducing kernel for the space $F^2(\mathbb{C}^n)$. Namely,

$$f(z) = \int_{\mathbb{C}^n} f(w) K(z, w) e^{-|w|^2} \, dV(w), \qquad z \in \mathbb{C}^n \qquad (3.4)$$

<span style="float:right">`repro`</span>

for $f \in F^2(\mathbb{C}^n)$; see, for example, [11, Proposition 2.2] for one variable case.

We need some more notation. For $s$ real and $f \in H(\mathbb{C}^n)$, let $f_s^+$ be the tail part of the Taylor expansion of $f$ of degree bigger than $|s|$ and $f_s^- = f - f_s^+$. So, if (3.1) holds, then we have

$$f_s^+ = \sum_{k > |s|} f_k \quad \text{and} \quad f_s^- = \sum_{k \leq |s|} f_k. \qquad (3.5)$$

<span style="float:right">`taylor`</span>

For an integer $k \geq 0$, we denote by $e_k$ the $k$-th "truncated" exponential function given by

$$e_k(\lambda) = e^\lambda - \sum_{j=0}^k \frac{\lambda^j}{j!}, \qquad \lambda \in \mathbb{C}.$$

It is easy to check that

$$\frac{e_k(\lambda)}{\lambda^{k+1}} = \sum_{\ell=0}^\infty \frac{\lambda^\ell}{(k+1+\ell)!} = \frac{1}{k!} \int_0^1 (1-t)^k e^{t\lambda} dt, \qquad (3.6)$$

<span style="float:right">`ekl`</span>

which immediately yields a useful inequality

$$|e_k(\lambda)| \leq \left( \frac{|\lambda|}{\operatorname{Re}\lambda} \right)^{k+1} e_k(\operatorname{Re}\lambda) \qquad (3.7)$$

<span style="float:right">`compa`</span>

for $\lambda \in \mathbb{C}$. Also, we have

$$0 < \frac{e_k(x)}{x^{k+1}} \leq e^x \qquad (3.8)$$

<span style="float:right">`compb`</span>

for $x > 0$.

We now proceed to estimate the fractional derivatives of the Fock kernel. We begin with the integral representation for the fractional derivatives. In what follows $\partial_t := \frac{\partial}{\partial t}$.

<span style="float:left">`radialder`</span>
**Lemma 3.1.** *Let $s > 0$ and put $s = m + r$ where $m$ is a nonnegative integer and $0 \leq r < 1$. Then the following identities hold for $f \in H(\mathbb{C}^n)$ and $z \in \mathbb{C}^n$:*

$$\mathcal{D}^s f(z) = \begin{cases} m! f(0) + \displaystyle\int_0^1 \partial_t^{m+1}[t^m f(tz)] \, dt & \text{if } n = 1 \text{ and } r = 0 \\[2ex] \dfrac{1}{\Gamma(1-r)} \displaystyle\int_0^1 \dfrac{\partial_t^{m+1}[t^{n+s-1} f(tz)]}{(1-t)^r} \, dt & \text{otherwise} \end{cases}$$

*and*

$$\mathcal{D}^{-s} f(z) = \frac{1}{\Gamma(s)} \int_0^1 t^{n-s-1} (1-t)^{s-1} f_s^+(tz) \, dt.$$

*Proof.* We provide a proof for $\mathcal{D}^s$; the proof for $\mathcal{D}^{-s}$ is simpler. Using the homogeneous expansion of an entire function, we only need to prove the integral representation for homogeneous polynomials. So, assume that $f$ is a homogeneous polynomial of degree $k$ in the rest of the proof.

Fix $z \in \mathbb{C}^n$. When $n + r - 1 + k > 0$, note

$$\partial_t^{m+1}[t^{n+s-1}f(tz)] = \partial_t^{m+1}[t^{n+s-1+k}]f(z)$$
$$= \frac{\Gamma(n+s+k)}{\Gamma(n+r-1+k)}t^{n+r+k-2}f(z).$$

So, multiplying both sides by $(1-t)^{-r}/\Gamma(1-r)$ and then integrating, we obtain

$$\frac{1}{\Gamma(1-r)}\int_0^1 \frac{\partial_t^{m+1}[t^{n+s-1}f(tz)]}{(1-t)^r}\,dt$$
$$= f(z)\frac{\Gamma(n+s+k)}{\Gamma(n+r-1+k)\Gamma(1-r)}\int_0^1 (1-t)^{-r}t^{n+r+k-2}\,dt$$
$$= \frac{\Gamma(n+s+k)}{\Gamma(n+k)}f(z).$$

This completes the proof for the case when $n \geq 2$ or $0 < r < 1$, because $n+r-1+k > 0$ for all $k \geq 0$. The case when $n = 1$ and $r = 0$ is treated similarly, because the above integral representation remains valid for all $k \geq 1$ and $\mathcal{D}^m 1 = m!$. The proof is complete. $\qquad\square$

Given $\delta > 0$, put

$$A_\delta(z) := \{w \in \mathbb{C}^n : |\theta(z,w)| < \delta\}$$

for $z \in \mathbb{C}^n$ where $\theta(z,w)$ is the angle between $z$ and $w$ identified as real vectors in $\mathbb{R}^{2n}$ so that $\mathrm{Re}\,(z \cdot \overline{w}) = |z||w|\cos\theta(z,w)$. Also, given $\epsilon > 0$, put

$$\Lambda_{\epsilon,\delta}(z,w) := e^{\mathrm{Re}\,(z\cdot\overline{w})}\chi_{A_\delta(z)}(w) + e^{\epsilon|z||w|} \qquad (3.9) \quad \boxed{\texttt{sete}}$$

for $z, w \in \mathbb{C}^n$ where $\chi$ denotes the characteristic function of the set specified in the subscript. With these notation we have the following pointwise size estimate for the fractional derivatives of the Fock kernel.

$\boxed{\texttt{dsker}}$ **Proposition 3.2.** *Given $0 < \epsilon < 1$ and $s$ real, there are positive constants $C = C(s,\epsilon) > 0$ and $\delta = \delta(\epsilon) > 0$ such that*

$$|\mathcal{D}^s K_w(z)| \leq C \times \begin{cases} (1+|z\cdot\overline{w}|)^s\Lambda_{\epsilon,\delta}(z,w) & \text{if} \quad s > 0 \\ (1+|z||w|)^s\Lambda_{\epsilon,\delta}(z,w) & \text{if} \quad s < 0 \end{cases}$$

*for $z, w \in \mathbb{C}^n$.*

*Proof.* Fix $0 < \epsilon < 1$ and $s > 0$. Put $s = m + r$ where $m$ is a nonnegative integer and $0 \leq r < 1$. Given $z, w \in \mathbb{C}^n$, put $\lambda = z \cdot \overline{w}$ and $x = \mathrm{Re}\,\lambda$ for short.

First, we estimate $|\mathcal{D}^s K_w(z)|$. Our proof is based on the integral representation given in Lemma 3.1. We provide details only for the case when $n \geq 2$ or $0 < r < 1$; the remaining case is treated similarly. Since $\partial_t^{m+1}[t^{n+s-1}e^{t\lambda}]$ is equal to $e^{t\lambda}$ times a linear combination of $t^{n+j+r-2}\lambda^j$ with $j = 0, 1, \ldots, m+1$, we have by Lemma 3.1

$$|\mathcal{D}^s K_w(z)| \lesssim \sum_{j=0}^{m+1}(1+|\lambda|)^j \int_0^1 e^{tx}(1-t)^{-r}t^{n+j+r-2}\, dt. \qquad (3.10) \quad \boxed{\text{dskw}}$$

So, in case $x \le 1$, we have by (3.10)

$$\begin{aligned}
|\mathcal{D}^s K_w(z)| &\lesssim (1+|\lambda|)^{m+1} \\
&= (1+|\lambda|)^s e^{\epsilon|z||w|} \cdot \frac{(1+|\lambda|)^{1-r}}{e^{\epsilon|z||w|}} \\
&\lesssim (1+|\lambda|)^s e^{\epsilon|z||w|},
\end{aligned}$$

which implies the asserted estimate.

Now, assume $x > 1$. The first term of the sum in (3.10) is easily seen to be dominated by some constant times $e^x$. Meanwhile, the other terms are all dominated by some constant times

$$(1+|\lambda|^{m+1}) \int_0^1 e^{tx}(1-t)^{-r}\, dt = (1+|\lambda|^{m+1})\frac{e^x}{x^{1-r}} \int_0^x e^{-t}t^{-r}\, dt.$$

Note that the integral in the right-hand side of the above is bounded by $\int_0^\infty e^{-t}t^{-r}\, dt$, which is finite. Overall, we see from (3.10) that

$$|\mathcal{D}^s K_w(z)| \lesssim \frac{(1+|\lambda|)^{m+1}}{(1+x)^{1-r}}e^x = \left(\frac{1+|\lambda|}{1+x}\right)^{m+1}(1+x)^s e^x \qquad (3.11) \quad \boxed{\text{xx}}$$

for $x > 1$.

Now, choose $\delta \in (0, \pi/2)$ such that $2\cos\delta = \epsilon$. It is easily seen from (3.11) that the required estimate holds when $w \in A_\delta(z)$, because $x \approx |z||w| \approx |\lambda|$ for such $w$. So, assume $w \notin A_\delta(z)$. Note $x \le |z||w|\cos\delta$ for such $w$. We thus have by our choice of $\delta$

$$\begin{aligned}
\frac{(1+|\lambda|)^{m+1}}{(1+x)^{1-r}}e^x &\le \frac{(1+|\lambda|)^{m+1}}{(1+x)^{1-r}}e^{\epsilon|z||w|/2} \\
&\le (1+|\lambda|)^s e^{\epsilon|z||w|} \cdot \frac{(1+|\lambda|)^{1-r}}{e^{\epsilon|z||w|/2}} \qquad (3.12) \quad \boxed{\text{ds}} \\
&\lesssim (1+|\lambda|)^s e^{\epsilon|z||w|}.
\end{aligned}$$

This, together with (3.11), yields the asserted estimate for $x > 1$. This completes the proof for the case $s > 0$.

Next, we estimate $|\mathcal{D}^{-s}K_w(z)|$. In this case we have by Lemma 3.1

$$|\mathcal{D}^{-s}K_w(z)| \le \frac{1}{\Gamma(s)} \int_0^1 t^{n-s-1}(1-t)^{s-1}|e_m(t\lambda)|\, dt.$$

Thus we obtain by (3.7)

$$|\mathcal{D}^{-s}K_w(z)| \le \left(\frac{|\lambda|}{x}\right)^{m+1}\frac{1}{\Gamma(s)} \int_0^1 t^{n-s-1}(1-t)^{s-1}e_m(tx)\, dt. \qquad (3.13) \quad \boxed{\text{isint}}$$

Note from (3.6) that $|e_m(tx)| \lesssim (t|x|)^{m+1}$ when $x$ stays bounded above. So, in case $x \leq 2$, we have

$$|\mathcal{D}^{-s}K_w(z)| \lesssim |\lambda|^{m+1} \lesssim (1+|z||w|)^{-s}e^{\epsilon|z||w|},$$

which implies the asserted estimate.

Now, assume $x > 2$. Note $e_m(tx) > 0$. Denoting by $I$ the integral in the right-hand side of (3.13), we claim

$$I \leq Cx^{-s}e^x, \qquad x > 2 \qquad (3.14) \quad \boxed{\texttt{claim}}$$

for some constant $C > 0$ independent of $x$. In order to prove this claim, we write $I$ as the sum of three pieces $I_1$, $I_2$ and $I_3$ defined by

$$I_1 = \int_0^{1/x}, \qquad I_2 = \int_{1/x}^{1/2}, \qquad I_3 = \int_{1/2}^1$$

and show that each of these pieces satisfies the desired estimate. For the first integral, we note from (3.8) that $e_m(tx) < (tx)^{m+1}e$ for $0 < t < 1/x$. Thus we have

$$I_1 \lesssim \int_0^{1/x} (tx)^{m+1}t^{n-s-1}(1-t)^{s-1}\,dt \lesssim x^{m+1}.$$

For the second integral, we note from the definition of $e_m$ that $e_m(tx) < e^{tx} < e^{\frac{x}{2}}$ for $0 < t < 1/2$. Thus we have

$$I_2 \lesssim e^{x/2} \int_{1/x}^{1/2} t^{-s}\,dt \lesssim x^s e^{\frac{x}{2}}.$$

For the third integral, we have

$$I_3 \lesssim \int_0^1 (1-t)^{s-1}e^{tx}\,dt = x^{-s}e^x \int_0^x \tau^{s-1}e^{-\tau}\,d\tau \lesssim x^{-s}e^x.$$

Combining these observations together, we obtain $I \lesssim x^{m+1} + x^s e^{\frac{x}{2}} + x^{-s}e^x$ for all $x$, which implies (3.14).

Now, having (3.14), we see from (3.13) that

$$|\mathcal{D}^{-s}K_w(z)| \lesssim \left(\frac{|\lambda|}{x}\right)^{m+1} x^{-s}e^x, \qquad x > 2.$$

Form this it is easily verified that the asserted estimate for $|\mathcal{D}^{-s}K_w(z)|$ holds for $w \in A_\delta(z)$, as in the case of $|\mathcal{D}^s K_w(z)|$. Also, for $w \notin A_\delta(z)$, we have as in the argument of (3.12)

$$\left(\frac{|\lambda|}{x}\right)^{m+1} x^{-s}e^x \lesssim |\lambda|^{m+1}e^{\epsilon|z||w|/2} \lesssim (1+|z||w|)^{-s}e^{\epsilon|z||w|}, \qquad (3.15) \quad \boxed{\texttt{xbig}}$$

which implies the asserted estimate. The proof is complete. $\qquad \square$

We now estimate the fractional integrals of the Fock kernel. For that purpose we need a couple of lemmas. First, we observe that the fractional integrals also admit integral representations, as in the case of the fractional derivatives.

`radialint` **Lemma 3.3.** *Let $s > 0$ and put $s = m + r$ where $m$ is a nonnegative integer and $0 \leq r < 1$. Then the following identities hold for $f \in H(\mathbb{C}^n)$ and $z \in \mathbb{C}^n$:*

$$\mathcal{I}^s f(z) = \frac{1}{\Gamma(s)} \int_0^1 t^{n-1}(1-t)^{s-1} f(tz)\, dt$$

*and*

$$\mathcal{I}^{-s} f(z) = \frac{1}{\Gamma(1-r)} \int_0^1 \frac{t^s \partial_t^{m+1}[t^{n-r} f_s^+(tz)]}{(1-t)^r}\, dt.$$

*Proof.* We prove the second part; the proof for the first part is simpler. Let $f \in H(\mathbb{C}^n)$. As in the proof of Lemma 3.1, we may assume that $f$ is a homogeneous polynomial, say, of degree $k$. We may further assume $k \geq m + 1$ so that $f_s^+ = f$ to avoid triviality. Now, for $z \in \mathbb{C}^n$, since

$$\partial_t^{m+1}[t^{n-r} f(tz)] = \partial_t^{m+1}[t^{n+k-r}] f(z)$$
$$= \frac{\Gamma(n+k+1-r)}{\Gamma(n+k-s)} t^{n+k-s-1} f(z),$$

we obtain

$$\frac{1}{\Gamma(1-r)} \int_0^1 \frac{t^s \partial_t^{m+1}[t^{n-r} f_s^+(tz)]}{(1-t)^r}\, dt$$
$$= \frac{f(z)\Gamma(n+k+1-r)}{\Gamma(n+k-s)\Gamma(1-r)} \int_0^1 (1-t)^{-r} t^{n+k-1}\, dt$$
$$= \frac{\Gamma(n+k)}{\Gamma(n+k-s)} f(z),$$

as required. The proof is complete. $\qquad\square$

Next, we need the following information on the derivatives of the truncated exponential functions.

`e-dif` **Lemma 3.4.** *Given $a$ real and an integer $m \geq 1$, there is a constant $C = C(a, m) > 0$ such that*

$$\left| \partial_t^{m+1}[t^a e_m(t\lambda)] \right| \leq C\, t^a |\lambda|^{m+1} e^{t\operatorname{Re}\lambda}$$

*for $t > 0$ and $\lambda \in \mathbb{C}$ with $\operatorname{Re} \lambda > 0$.*

*Proof.* Fix a real number $a$ and an integer $m \geq 1$. Let $\lambda \in \mathbb{C}$ with $x := \operatorname{Re} \lambda > 0$. Since $e_k' = e_{k-1}$ for integers $k \geq 0$ where $e_{-1}$ is the original exponential function, we see that $\partial_t^{m+1}[t^a e_m(t\lambda)]$ is equal to a linear combination of $t^{a-j}\lambda^{m+1-j} e_{j-1}(t\lambda)$ with $j = 0, 1, \ldots, m+1$. We thus have by (3.7) and (3.8)

$$\left| \partial_t^{m+1}[t^a e_m(t\lambda)] \right| \lesssim t^a |\lambda|^{m+1} \left\{ e^{tx} + \sum_{j=1}^{m+1} \frac{|e_{j-1}(tx)|}{|tx|^j} \right\}$$
$$\leq (m+3) t^a |\lambda|^{m+1} e^{tx},$$

as required. The proof is complete. $\qquad\square$

We are now ready to prove the following pointwise size estimate for the fractional integrals of the Fock kernel.

**Proposition 3.5.** *Given* $0 < \epsilon < 1$ *and* $s$ *real, there are constants* $C = C(s, \epsilon) > 0$ *and* $\delta = \delta(\epsilon) > 0$ *such that*

$$|\mathcal{I}^s K_w(z)| \leq C \times \begin{cases} (1 + |z||w|)^{-s} \Lambda_{\epsilon, \delta}(z, w) & \text{if } s > 0 \\ |z \cdot \overline{w}|^{-s} \Lambda_{\epsilon, \delta}(z, w) & \text{if } s < 0 \end{cases}$$

*for* $z, w \in \mathbb{C}^n$.

*Proof.* Fix $0 < \epsilon < 1$ and $s > 0$. Put $s = m + r$ where $m$ is a nonnegative integer and $0 \leq r < 1$. Given $z, w \in \mathbb{C}^n$, we continue using the notation introduced in the proof of Proposition 3.2. So, $\lambda = z \cdot \overline{w}$ and $x = \operatorname{Re} \lambda$. Also, $\delta \in (0, \pi/2)$ is chosen so that $2 \cos \delta = \epsilon$.

First, we estimate $|\mathcal{I}^s K_w(z)|$. We have by Lemma 3.3

$$|\mathcal{I}^s K_w(z)| \lesssim \int_0^1 (1-t)^{s-1} e^{tx} \, dt.$$

Note that the right-hand side of the above stays bounded for $x \leq 1$. Meanwhile, we have

$$\int_0^1 (1-t)^{s-1} e^{tx} \, dt \lesssim x^{-s} e^x.$$

for $x > 1$. Now, slightly modifying the argument for the estimate of $|\mathcal{D}^{-s} K_w(z)|$ in the proof of Proposition 3.2, we see that the asserted estimate holds.

Next, we estimate $|\mathcal{I}^{-s} K_w(z)|$. Note $(K_w)_s^+(tz) = e_m(t\lambda)$. Thus we have by Lemmas 3.3 and 3.4

$$|\mathcal{I}^{-s} K_w(z)| \lesssim |\lambda|^{m+1} \int_0^1 (1-t)^{-r} t^{m+n} e^{tx} \, dt.$$

Thus we have

$$|\mathcal{I}^{-s} K_w(z)| \lesssim |\lambda|^{m+1}$$

for $x \leq 1$. Meanwhile, for $x > 1$, we have

$$\begin{aligned} |\mathcal{I}^{-s} K_w(z)| &\lesssim |\lambda|^{m+1} \int_0^1 (1-t)^{-r} e^{tx} \, dt \\ &\lesssim |\lambda|^{m+1} x^{r-1} e^x \\ &= \left( \frac{|\lambda|}{x} \right)^{m+1} x^s e^x. \end{aligned}$$

Thus, slightly modifying the argument for the estimate of $|\mathcal{D}^s K_w(z)|$ in the proof of Proposition 3.2, we see that the asserted estimate holds. The proof is complete. $\square$

Having seen Propositions 3.2 and 3.5, we now turn to the $L^p$-integral estimates, with respect to weighted Gaussian measures, for the functions $\Lambda_{\epsilon,\delta}$. Note

$$|\Lambda_{\epsilon,\delta}(z,w)| \leq e^{\operatorname{Re}(z\cdot\overline{w})} + e^{\epsilon|z||w|}, \qquad z,w \in \mathbb{C}^n \qquad (3.16)$$

for any $\delta, \epsilon > 0$. So, we consider $L^p$-integrals of each term in the right-hand-side of the above separately. First, for the first term of (3.16), we have the following integral estimate.

**Lemma 3.6.** *Given $0 < p, a < \infty$ and $\alpha$ real, there is a constant $C = C(p,a,\alpha) > 0$ such that*

$$\int_{\mathbb{C}^n} e^{p\operatorname{Re}(z\cdot\overline{w})-a|w|^2}\, dV_\alpha(w) \leq C\frac{e^{\frac{p^2}{4a}|z|^2}}{(1+|z|)^\alpha}$$

*for $z \in \mathbb{C}^n$.*

*Proof.* Let $0 < p, a < \infty$ and $\alpha$ be a real number. Given $z,w \in \mathbb{C}^n$, note

$$p\operatorname{Re}(z\cdot\overline{w}) - a|w|^2 = \frac{p^2}{4a}|z|^2 - \left|\frac{p}{2\sqrt{a}}z - \sqrt{a}w\right|^2.$$

Also, note by (2.1)

$$\begin{aligned}
\frac{1}{(1+|w|)^\alpha} &\approx \frac{1}{(1+\frac{p}{2\sqrt{a}}|w|)^\alpha} \\
&\leq \frac{1}{\left(1+\frac{p}{2\sqrt{a}}|z|\right)^\alpha}\left(1+\left|\frac{p}{2\sqrt{a}}z-\sqrt{a}w\right|\right)^{|\alpha|} \\
&\approx \frac{1}{(1+|z|)^\alpha}\left(1+\left|\frac{p}{2\sqrt{a}}z-\sqrt{a}w\right|\right)^{|\alpha|}.
\end{aligned}$$

It follows that the integral under consideration is dominated by some constant times

$$\frac{e^{\frac{p^2}{4a}|z|^2}}{(1+|z|)^\alpha}\int_{\mathbb{C}^n} e^{-|\xi|^2}(1+|\xi|)^{|\alpha|}\, dV(\xi).$$

Now, since the integral above is finite, we conclude the lemma. $\qquad\square$

Next, for the second term of (3.16), we have the following integral estimate.

**Lemma 3.7.** *Given $0 < p, a, \epsilon < \infty$ and $\alpha$ real, there is a constant $C = C(p,a,\epsilon,\alpha) > 0$ such that*

$$\int_{\mathbb{C}^n} e^{p\epsilon|z||w|-a|w|^2}\, dV_\alpha(w) \leq Ce^{\frac{p^2\epsilon^2}{4a}|z|^2} \times \begin{cases} 1+|z|^{2n-\alpha} & \text{if } \alpha \neq 2n \\ 1+\log(1+|z|) & \text{if } \alpha = 2n \end{cases}$$

*for $z \in \mathbb{C}^n$.*

*Proof.* Let $0 < p, a, \epsilon < \infty$ and $\alpha$ be a real number. Given $z \in \mathbb{C}^n$, denote by $I(z)$ the integral under consideration. Note

$$I(z) = e^{\frac{p^2\epsilon^2}{4a}|z|^2} \int_{\mathbb{C}^n} e^{-a\left(\frac{p\epsilon}{2a}|z|-|w|\right)^2} \, dV_\alpha(w)$$

$$=: e^{\frac{p^2\epsilon^2}{4a}|z|^2} [I_1(z) + I_2(z)]$$

where

$$I_1(z) = \int_{|w| \le p\epsilon|z|/a} \qquad \text{and} \qquad I_2(z) = \int_{|w| > p\epsilon|z|/a}.$$

Since

$$I_1(z) \le \int_{|w| \le p\epsilon|z|/a} dV_\alpha(w),$$

an integration in polar coordinates yields

$$I_1(z) \lesssim \begin{cases} 1 + |z|^{2n-\alpha} & \text{if} \quad \alpha \ne 2n \\ \log(1 + |z|) & \text{if} \quad \alpha = 2n. \end{cases}$$

Meanwhile, since $|w| - p\epsilon|z|/2a \ge |w|/2$ for $|w| > p\epsilon|z|/a$, we have for the second integral

$$I_2(z) \le \int_{\mathbb{C}^n} e^{-\frac{a|w|^2}{4}} \, dV_\alpha(w) \approx 1.$$

Combining these observations together, we conclude the lemma. $\qquad\square$

Now, as an immediate consequence of (3.16), Lemmas 3.6 and 3.7, we have the next estimate for the $L^p$-integrals of $\Lambda_{\epsilon,\delta}$ against weighted Gaussian measures, when restricted to $0 < \epsilon < 1$. The next estimate turns out to be enough for our purpose, although it is derived from the very rough inequality (3.16).

ebound **Proposition 3.8.** *Given $0 < p, a < \infty$, $0 < \epsilon < 1$ and $\alpha$ real, there is a constant $C = C(p, a, \epsilon, \alpha) > 0$ such that*

$$\int_{\mathbb{C}^n} |\Lambda_{\epsilon,\delta}(z,w)|^p e^{-a|w|^2} \, dV_\alpha(w) \le C \frac{e^{\frac{p^2}{4a}|z|^2}}{(1 + |z|)^\alpha}$$

*for $\delta > 0$ and $z \in \mathbb{C}^n$.*

We now proceed to investigating how the fractional differentiation/integration acts on the weighted Fock spaces.

To handle the case $1 \le p < \infty$ and for other purpose later, we introduce an auxiliary class of integral operators. Fix $0 < \epsilon < 1$ and $\delta > 0$. Given $s$ real, we consider an integral operator $L_s = L_{s,\epsilon,\delta}$ defined by

$$L_s\psi(z) := \int_{\mathbb{C}^n} \psi(w) \left(\frac{1 + |z|}{1 + |w|}\right)^s \Lambda_{\epsilon,\delta}(z,w) e^{-|w|^2} \, dV(w), \qquad z \in \mathbb{C}^n$$

for $\psi$ which makes the above integral well-defined.

`lbound` **Proposition 3.9.** *Given $s$ real, the operator $L_s$ is bounded on $L_\alpha^p$ for any $1 \le p \le \infty$ and $\alpha$ real.*

*Proof.* Fix a real number $\alpha$. We first consider the case $s = 0$. Put $\Lambda := \Lambda_{\epsilon,\delta}$ for short. By Fubini's theorem we have

$$
\begin{aligned}
\|L_0 \psi\|_{L_\alpha^1} &= \int_{\mathbb{C}^n} \left| \int_{\mathbb{C}^n} \Lambda(z,w)\psi(w)e^{-|w|^2}\, dV(w) \right| e^{-\frac{1}{2}|z|^2}\, dV_\alpha(z) \\
&\le \int_{\mathbb{C}^n} |\psi(w)| e^{-|w|^2} \left\{ \int_{\mathbb{C}^n} \Lambda(z,w) e^{-\frac{1}{2}|z|^2}\, dV_\alpha(z) \right\} dV(w)
\end{aligned}
$$

for $\psi \in L_\alpha^1$. Since the inner integral of the above is dominated by some constant times $e^{\frac{1}{2}|w|^2}(1+|w|)^{-\alpha}$ by Proposition 3.8, we see that $L_0$ is bounded on $L_\alpha^1$.

Next, we have again by Proposition 3.8

$$
\begin{aligned}
|L_0 \psi(z)| &\lesssim \|\psi\|_{L_\alpha^\infty} \int_{\mathbb{C}^n} \Lambda(z,w) e^{-\frac{1}{2}|w|^2}\, dV_{-\alpha}(w), \qquad z \in \mathbb{C}^n \\
&\lesssim e^{\frac{1}{2}|z|^2}(1+|z|)^\alpha \|\psi\|_{L_\alpha^\infty}
\end{aligned}
$$

for $\psi \in L_\alpha^\infty$. So, $L_0$ is bounded on $L_\alpha^\infty$. In particular, $L_0$ is bounded on $L_0^\infty$. Thus it follows from the Stein interpolation theorem (see [3, Theorem 3.6]) that $L_0$ is bounded on $L_\alpha^p$ for any $1 \le p < \infty$. This completes the proof for $s = 0$.

Now, we consider general $s$. Note

$$
\begin{aligned}
\frac{L_s \psi(z)}{(1+|z|)^s} &= \int_{\mathbb{C}^n} \frac{\psi(w)}{(1+|w|)^s} \Lambda(z,w) e^{-|w|^2}\, dV(w) \\
&= L_0 \left[ \frac{\psi(w)}{(1+|w|)^s} \right](z).
\end{aligned}
$$

Thus, for $1 \le p < \infty$, we see that $L_s$ is bounded on $L_\alpha^p$ by the boundedness of $L_0$ on $L_{\alpha-ps}^p$. Also, we see that $L_s$ is bounded on $L_\alpha^\infty$ by the boundedness of $L_0$ on $L_{\alpha-s}^p$. The proof is complete. $\qquad \square$

The following Jensen-type inequality is needed to handle the case $0 < p \le 1$.

`smmvplem` **Proposition 3.10.** *Given $0 < p \le 1$, $a > 0$ and $\alpha$ real, there is a constant $C = C(p,a,\alpha) > 0$ such that*

$$
\left\{ \int_{\mathbb{C}^n} |f(z)| e^{-a|z|^2}\, dV_\alpha(z) \right\}^p \le C \int_{\mathbb{C}^n} \left| f(z) e^{-a|z|^2} \right|^p\, dV_{p\alpha}(z) \qquad (3.17)
$$  `mvpsmall`

*for $f \in H(\mathbb{C}^n)$.*

*Proof.* Let $0 < p \le 1$, $a > 0$ and $\alpha$ be a real number. Let $f \in H(\mathbb{C}^n)$. By Lemma 2.1 there is a constant $C = C(p,a,\alpha) > 0$

$$
\frac{|f(z)|}{e^{a|z|^2}(1+|z|)^\alpha} \le C \left\{ \int_{\mathbb{C}^n} \left| f(w) e^{-a|w|^2} \right|^p\, dV_{\alpha p}(w) \right\}^{1/p}
$$

and hence

$$\frac{|f(z)|}{e^{a|z|^2}(1+|z|)^\alpha} = \left|\frac{f(z)}{e^{a|z|^2}(1+|z|)^\alpha}\right|^p \left|\frac{f(z)}{e^{a|z|^2}(1+|z|)^\alpha}\right|^{1-p}$$

$$\lesssim \left|\frac{f(z)}{e^{a|z|^2}(1+|z|)^\alpha}\right|^p \left\{\int_{\mathbb{C}^n} \left|f(w)e^{-a|w|^2}\right|^p dV_{\alpha p}(w)\right\}^{(1-p)/p}$$

for $z \in \mathbb{C}^n$. Now, integrating both sides of the above against the measure $dV(z)$, we conclude the proposition. $\square$

Given $0 < p < \infty$ and $\alpha$ real, it is not hard to see via the subharmonicity and the maximum modulus theorem that $\sup_{|w| \leq 1} |f(w)|$ is dominated by some constant times $\|f\|_{F_\alpha^p}$ for any $f \in F_\alpha^p$. This property extends to arbitrary fractional derivatives as in the next lemma.

**Lemma 3.11.** *Given $0 < p < \infty$ and $s, \alpha$ real, there is a constant $C = C(p, s, \alpha) > 0$ such that*

$$\sup_{|z| \leq 1} |\mathcal{D}^s f(z)| + \sup_{|z| \leq 1} |\mathcal{I}^s f(z)| \leq C\|f\|_{F_\alpha^p}$$

*for $f \in F_\alpha^p$.*

*Proof.* We provide a proof only for $\mathcal{D}^s$; the proof for $\mathcal{I}^s$ is similar. Let $0 < p < \infty$ and $s, \alpha$ be real numbers. By Proposition 2.3 it is sufficient to consider only holomorphic polynomials. So, fix an arbitrary holomorphic polynomial $f$. Applying $\mathcal{D}^s$ to (3.4), we have

$$\mathcal{D}^s f(z) = \int_{\mathbb{C}^n} f(w)\mathcal{D}^s K_w(z)e^{-|w|^2} dV(w)$$

and thus

$$|\mathcal{D}^s f(z)| \leq \int_{\mathbb{C}^n} |f(w)\mathcal{D}^s K_w(z)|e^{-|w|^2} dV(w) \tag{3.18}$$

for $z \in \mathbb{C}^n$.

We now consider the cases $0 < p < 1$ and $1 \leq p < \infty$ separately. Assume $0 < p < 1$. Applying Proposition 3.10 to the holomorphic function $f(w)\overline{\mathcal{D}^s K_w(z)}$ with $z$ fixed, we obtain from (3.18)

$$|\mathcal{D}^s f(z)|^p \lesssim \int_{\mathbb{C}^n} |f(w)\mathcal{D}^s K_w(z)|^p e^{-p|w|^2} dV(w) =: I(z) \tag{3.19}$$

for $z \in \mathbb{C}^n$. Note

$$|\mathcal{D}^s K_w(z)| \lesssim (1 + |z||w|)^s e^{|z||w|}, \qquad z, w \in \mathbb{C}^n \tag{3.20}$$

by Proposition 3.2 and (3.16). Thus, for $|z| \leq 1$, we have

$$I(z) \lesssim \int_{\mathbb{C}^n} |f(w)|^p (1 + |w|)^{ps} e^{-p|w|^2 + p|w|} dV(w)$$

$$= \int_{\mathbb{C}^n} \left|f(w)e^{-\frac{|w|^2}{2}}\right|^p (1 + |w|)^{\alpha + ps} e^{-\frac{p|w|^2}{2} + p|w|} dV_\alpha(w)$$

$$\lesssim \|f\|_{L_\alpha^p}^p,$$

as desired.

Next, assume $1 \leq p < \infty$. We have by (3.18) and (3.20)

$$|\mathcal{D}^s f(z)| \leq \int_{\mathbb{C}^n} \left| f(w) e^{-\frac{|w|^2}{2}} \right| (1 + |w|)^s e^{-\frac{|w|^2}{2} + |w|} \, dV(w)$$

for $|z| \leq 1$. Thus, applying Jensen's inequality with respect to the finite measure $d\mu(w) := (1 + |w|)^s e^{-\frac{|w|^2}{2} + |w|} \, dV(w)$, we obtain

$$|\mathcal{D}^s f(z)| \lesssim \int_{\mathbb{C}^n} \left| f(w) e^{-\frac{|w|^2}{2}} \right|^p \, d\mu(w) \lesssim \|f\|_{L^p_\alpha}^p,$$

for $|z| \leq 1$. This completes the proof. $\qquad\square$

**Lemma 3.12.** *Given $\alpha$ real and $a, b > 0$, there is a constant $C = C(\alpha, a, b) > 0$ such that*

$$\int_0^1 t^{a-1}(1-t)^{b-1} e^{|tz|^2/2}(1 + |tz|)^\alpha \, dt \leq C e^{|z|^2/2}(1 + |z|)^{\alpha - 2b}$$

*for $z \in \mathbb{C}^n$.*

*Proof.* Denote by $I(z)$ the integral in question. Since $I(z)$ stays bounded for $|z| \leq 1$, we may assume $|z| \geq 1$.

Decompose $I(z)$ into two pieces

$$I(z) = \int_0^{1/2} + \int_{1/2}^1 .$$

The first integral is easily treated, because $I_1(z) \lesssim (1 + |z|)^\alpha e^{\frac{|z|^2}{8}}$ if $\alpha \geq 0$ and $I_1(z) \lesssim e^{\frac{|z|^2}{8}}$ if $\alpha < 0$. Since $1 + |tz| \approx 1 + |z|$ for $1/2 \leq t < 1$, we have

$$\int_{1/2}^1 \approx (1 + |z|)^\alpha \int_{1/2}^1 (1 - t^2)^{b-1} e^{|tz|^2/2} \, dt$$

$$\approx e^{|z|^2/2}(1 + |z|)^\alpha \int_0^{3/4} t^{b-1} e^{-t|z|^2/2} \, dt.$$

Meanwhile, we have

$$\int_0^{3/4} t^{b-1} e^{-t|z|^2/2} \, dt = |z|^{-2b} \int_0^{3|z|^2/4} x^{b-1} e^{-x/2} \, dx \lesssim (1 + |z|)^{-2b}.$$

Thus the required estimate holds. The proof is complete. $\qquad\square$

We are now ready to prove that each fractional differentiation/integration on a weighted Fock space amounts to increasing the weight as in the next theorem.

**Theorem 3.13.** *Let $s$ and $\alpha$ be real numbers. Then the operators*

$$\mathcal{D}^s, \mathcal{I}^{-s} : F^p_\alpha \to \begin{cases} F^p_{\alpha + 2sp} & \text{if } 0 < p < \infty \\ F^p_{\alpha + 2s} & \text{if } p = \infty \end{cases}$$

*are bounded.*

*Proof.* Fix real numbers $s$ and $\alpha$. We consider two cases $0 < p < \infty$ and $p = \infty$ separately.

The case $0 < p < \infty$: We provide a proof only for $\mathcal{D}^s$; the proof for $\mathcal{I}^{-s}$ is similar (with the help of Proposition 3.5 instead of Proposition 3.2 in this case). By Proposition 2.3 and Lemma 3.11, it suffices to produce a constant $C = C(p, s, \alpha) > 0$ such that

$$J := \int_{|z| \geq 1} \left| \mathcal{D}^s f(z) e^{-\frac{1}{2}|z|^2} \right|^p \, dV_{\alpha + 2ps}(z) \leq C \|f\|_{F_\alpha^p}^p \qquad (3.21)$$

for holomorphic polynomials $f$. So, fix an arbitrary holomorphic polynomial $f$ and let $\Lambda = \Lambda_{\epsilon,\delta}$ be the function provided by Proposition 3.2 with $0 < \epsilon < 1$ fixed.

We now consider the cases $0 < p \leq 1$ and $1 < p < \infty$ separately. Assume $0 < p \leq 1$. In this case we have (3.19). Let $I(z)$ be the integral defined in (3.19) and decompose

$$I(z) = \int_{|w| \leq 1} + \int_{|w| > 1} =: I_1(z) + I_2(z).$$

For the first term, we have by (3.20) and Lemma 3.11 (with $s = 0$)

$$I_1(z) \lesssim (1 + |z|)^{p|s|} e^{p|z|} \|f\|_{F_\alpha^p}^p.$$

Meanwhile, note by Proposition 3.2

$$I_2(z) \lesssim \int_{|w| > 1} |f(w)|^p \, (1 + |z||w|)^{sp} |\Lambda(z, w)|^p e^{-p|w|^2} \, dV(w).$$

Since $(1 + |z||w|) \approx (1 + |z|)(1 + |w|)$ for $|z| \geq 1$ and $|w| \geq 1$, the above yields

$$I_2(z) \lesssim Q_{p,s} f(z), \qquad |z| \geq 1$$

where

$$Q_{p,s} f(z) := (1 + |z|)^{ps} \int_{\mathbb{C}^n} \left| f(w) e^{-|w|^2} \right|^p |\Lambda(z, w)|^p \, dV_{-ps}(w).$$

Combining these observations, we have so far

$$|\mathcal{D}^s f(z)|^p \lesssim (1 + |z|)^{p|s|} e^{p|z|} \|f\|_{F_\alpha^p}^p + Q_{p,s} f(z) \qquad (3.22)$$

for $|z| \geq 1$. Note

$$\int_{|z| \geq 1} (1 + |z|)^{p|s|} e^{-\frac{p|z|^2}{2} + p|z|} \, dV_{\alpha + 2ps}(z) < \infty,$$

which, together with (3.22), yields

$$J \lesssim \|f\|_{F_\alpha^p}^p + \int_{|z| \geq 1} Q_{p,s} f(z) e^{-\frac{p}{2}|z|^2} \, dV_{\alpha + 2ps}(z).$$

Note that the last integral is equal to

$$\int_{\mathbb{C}^n} \left| f(w) e^{-|w|^2} \right|^p \left\{ \int_{|z| \geq 1} |\Lambda(z, w)|^p e^{-\frac{p}{2}|z|^2} \, dV_{\alpha + ps}(z) \right\} dV_{-ps}(w),$$

which, in turn, is dominated by some constant times

$$\int_{\mathbb{C}^n} \left| f(w)e^{-|w|^2} \right|^p \frac{e^{\frac{p}{2}|w|^2}}{(1+|w|)^{\alpha+ps}} \, dV_{-ps}(w) = \|f\|_{F_\alpha^p}^p$$

by Proposition 3.8. So, (3.21) holds for $0 < p \le 1$.

Now, assume $1 < p < \infty$. Proceeding as in the case of $p = 1$ with the help of Lemma 3.11, we have

$$|\mathcal{D}^s f(z)| \lesssim (1+|z|)^{|s|} e^{|z|} \|f\|_{F_\alpha^p} + Q_{1,s} f(z), \qquad |z| \ge 1$$

and thus

$$J \lesssim \|f\|_{F_\alpha^p}^p + \|Q_{1,s} f\|_{L_{\alpha+2ps}^p}^p.$$

Meanwhile, since

$$\frac{Q_{1,s} f(z)}{(1+|z|)^{2s}} = \int_{\mathbb{C}^n} |f(w)| \left( \frac{1+|w|}{1+|z|} \right)^s \Lambda(z, w) e^{-|w|^2} \, dV(w),$$

we have $\|Q_{1,s} f\|_{L_{\alpha+2ps}^p} \lesssim \|f\|_{F_\alpha^p}$ by Proposition 3.9. So, (3.21) holds for $1 < p < \infty$. This completes the proof for $0 < p < \infty$.

The case $p = \infty$: In this case we assume $s > 0$ and provide proofs for $\mathcal{D}^s$ and $\overline{\mathcal{I}^{-s}}$; the proofs for other cases are simpler and the argument below can be easily modified. Write $s = m + r$ where $m$ is a nonnegative integer and $0 \le r < 1$. Let $f \in F_\alpha^\infty$.

First, we consider $\mathcal{D}^s$. Assume either $n \ge 2$ or $0 < r < 1$. Given $0 \le t \le 1$ and $z \in \mathbb{C}^n$, note

$$\begin{aligned}
\partial_t^{m+1}[t^{n+s-1} f(tz)] &= \sum_{j=0}^{m+1} c_{mj} t^{n+r-2+j} \partial_t^j [f(tz)] \\
&= t^{n+r-2} \sum_{j=0}^{m+1} j! c_{mj} \sum_{|\gamma|=j} \frac{(tz)^\gamma \partial^\gamma f(tz)}{\gamma!}
\end{aligned} \tag{3.23}$$

for some coefficients $c_{mj}$. Thus we have by Proposition 2.2

$$\begin{aligned}
|\partial_t^{m+1}[t^{n+s-1} f(tz)]| &\lesssim t^{n+r-2} \sum_{j=0}^{m+1} \sum_{|\gamma|=j} |\partial^\gamma f(tz)| |tz|^j \\
&\lesssim t^{n+r-2} e^{|tz|^2/2} (1+|tz|)^{\alpha+2(m+1)} \|f\|_{F_\alpha^\infty}.
\end{aligned}$$

So, we conclude by Lemmas 3.1 and 3.12 (with $a = n + r - 1$ and $b = 1 - r$)

$$\begin{aligned}
|\mathcal{D}^s f(z)| &\lesssim \|f\|_{F_\alpha^\infty} \int_0^1 t^{n+r-2}(1-t)^{-r} e^{|tz|^2/2}(1+|tz|)^{\alpha+2(m+1)} \, dt \\
&\lesssim e^{|z|^2/2}(1+|z|)^{\alpha+2(r-1+m+1)} \|f\|_{F_\alpha^\infty} \\
&= e^{|z|^2/2}(1+|z|)^{\alpha+2s} \|f\|_{F_\alpha^\infty}.
\end{aligned}$$

Now, assume $n = 1$ and $r = 0$. Choosing coefficients $a_{m\ell}$ such that $\frac{(k+m)!}{k!} = \sum_{\ell=0}^{m} a_{m\ell} k^\ell$ for all integers $k \geq 0$, we have

$$\mathcal{D}^m f(z) = \sum_{k=0}^{\infty} \left( \sum_{\ell=0}^{m} a_{m\ell} k^\ell \right) f_k(z) = \sum_{\ell=0}^{m} a_{m\ell} \left( z\frac{\partial}{\partial z} \right)^\ell f(z) \qquad (3.24)$$ `full-m`

and hence

$$|\mathcal{D}^m f(z)| \lesssim (1+|z|)^m \sum_{|\gamma| \leq m} |\partial^\gamma f(z)|$$

for all $z \in \mathbb{C}$. Thus we have the desired estimate by Proposition 2.2. This completes the proof for $\mathcal{D}^s$.

Now, we consider $\mathcal{I}^{-s}$. As in (3.23), we note

$$\partial_t^{m+1}[t^{n-r} f_s^+(tz)] = \sum_{j=0}^{m+1} j! c_{mj} t^{n-s-1} \sum_{|\gamma|=j} \frac{(tz)^\gamma \partial^\gamma f_s^+(tz)}{\gamma!}$$

and thus

$$|\partial_t^{m+1}[t^{n-r} f_s^+(tz)]| \lesssim t^{-s}(1+t|z|)^{m+1} \sum_{j=0}^{m+1} \sum_{|\gamma|=j} |\partial^\gamma f_s^+(tz)|$$

for $0 \leq t \leq 1$ and $z \in \mathbb{C}^n$. In order to estimate the size of the sum in the right-hand side of the above, we first note by Taylor's formula

$$\begin{aligned}
f_s^+(z) &= \frac{1}{m!} \int_0^1 (1-t)^m \partial_t^{m+1}[f(tz)]\, dt \\
&= (m+1) \sum_{|\nu|=m+1} \frac{z^\nu}{\nu!} \int_0^1 (1-t)^m \partial^\nu f(tz)\, dt. \qquad (3.25)
\end{aligned}$$ `taylor-1`

This, together with Proposition 2.2 and Lemma 3.12 (with $b = m+1$), yields

$$\begin{aligned}
|f_s^+(z)| &\lesssim |z|^{m+1} \sum_{|\nu|=m+1} \int_0^1 (1-t)^m |\partial^\nu f(tz)|\, dt \\
&\lesssim \|f\|_{F_\alpha^\infty} |z|^{m+1} \int_0^1 (1-t)^m e^{|tz|^2/2}(1+t|z|)^{\alpha+m+1}\, dt \\
&\lesssim |z|^{m+1} e^{|z|^2/2}(1+|z|)^{\alpha-m-1} \|f\|_{F_\alpha^\infty}.
\end{aligned}$$

In particular, we have $\|f_s^+\|_{F_\alpha^\infty} \lesssim \|f\|_{F_\alpha^\infty}$. Thus we deduce from Proposition 2.2

$$\sum_{j=0}^{m+1} \sum_{|\gamma|=j} |\partial^\gamma f_s^+(tz)| \lesssim e^{|tz|^2/2}(1+t|z|)^{\alpha+m+1} \|f\|_{F_\alpha^\infty}$$

so that

$$|\partial_t^{m+1}[t^{n-r} f_s^+(tz)]| \lesssim t^{-s} e^{|tz|^2/2}(1+t|z|)^{\alpha+2(m+1)} \|f\|_{F_\alpha^\infty}.$$

Now, as in the proof for $\mathcal{D}^s$, we obtain by Lemmas 3.3 and 3.12

$$|\mathcal{I}^{-s}f(z)| \lesssim e^{|z|^2/2}(1+|z|)^{\alpha+2s}\|f\|_{F_\alpha^\infty},$$

as required. The proof is complete. $\qquad\square$

Note that fractional derivatives and integrals of holomorphic polynomials are again holomorphic polynomials. Thus, as a consequence of Proposition 2.3 and Theorem 3.13, we also see that the operators

$$\mathcal{D}^s, \mathcal{I}^{-s} : F_\alpha^{\infty,0} \to F_{\alpha+2s}^{\infty,0}$$

are bounded for any $\alpha$, $s$ real.

## 4. Weighted Fock-Sobolev Spaces

In this section we introduce two types of weighted Fock-Sobolev spaces, one in terms of $\mathcal{R}^s$ and the other in terms of $\widetilde{\mathcal{R}}^s$. We first identify those spaces with the weighted Fock spaces. Then we describe explicitly the reproducing kernels.

Based on two notions of fractional differentiation/integration given in the previous section, we now introduce two different types of fractional radial differentiation/integration operators. For any $s$ real, we define the fractional radial differentiation/integration operators $\mathcal{R}^s$ and $\widetilde{\mathcal{R}}^s$ by

$$\mathcal{R}^s f(z) = \frac{1}{(1+|z|)^s}\mathcal{D}^s f(z) \tag{4.1}$$

and

$$\widetilde{\mathcal{R}}^s f(z) = \frac{1}{(1+|z|)^s}\mathcal{I}^{-s} f(z) \tag{4.2}$$

for $f \in H(\mathbb{C}^n)$. The weight factor $(1+|z|)^{-s}$ may look peculiar at first glance, but it plays an important normalization role in $\mathbb{C}^n$. In fact such a weight factor can be ignored on a bounded domain like the unit ball, as far as the growth behavior near boundary is concerned.

For $0 < p \le \infty$ and real numbers $\alpha$ and $s$, we define the weighted Fock-Sobolev space $F_{\alpha,s}^p$ to be the space of all $f \in H(\mathbb{C}^n)$ such that $\mathcal{R}^s f \in L_\alpha^p$ where $L_\alpha^p$ is the space introduced in the Introduction. We define the norm of $f \in F_{\alpha,s}^p$ by

$$\|f\|_{F_{\alpha,s}^p} := \begin{cases} \|\mathcal{R}^s f\|_{L_\alpha^p} & \text{if } s \ge 0 \\ \|\mathcal{R}^s f\|_{L_\alpha^p} + \|f_s^-\|_{F_\alpha^p} & \text{if } s < 0. \end{cases}$$

Similarly, the other type of weighted Fock-Sobolev $\widetilde{F}_{\alpha,s}^p$ is defined to be the space of all $f \in H(\mathbb{C}^n)$ such that $\widetilde{\mathcal{R}}^s f \in L_\alpha^p$ whose norm is given by

$$\|f\|_{\widetilde{F}_{\alpha,s}^p} := \begin{cases} \|\widetilde{\mathcal{R}}^s f\|_{L_\alpha^p} & \text{if } s \le 0 \\ \|\widetilde{\mathcal{R}}^s f\|_{L_\alpha^p} + \|f_s^-\|_{F_\alpha^p} & \text{if } s > 0. \end{cases}$$

Recall that $f_s^-$ is the Taylor polynomial of $f$ of degree less than or equal to $|s|$. In conjunction with these definitions we note for any parameters $p, \alpha$ and $s$

$$\|\mathcal{R}^s f\|_{L_\alpha^p} = \|\mathcal{D}^s f\|_{F_{\alpha+ps}^p} \quad \text{and} \quad \|\widetilde{\mathcal{R}}^s f\|_{L_\alpha^p} = \|\mathcal{I}^{-s} f\|_{F_{\alpha+ps}^p} \qquad (4.3)$$

`normrel`

for $f \in H(\mathbb{C}^n)$ with convention $F_{\alpha+ps}^p = F_{\alpha+s}^\infty$ for $p = \infty$.

`f-est`

**Lemma 4.1.** *Given $0 < p \le \infty$, real numbers $\alpha, \beta, s$ and a positive integer $m$, there is a constant $C = C(p, m, \alpha, \beta, s) > 0$ such that*

$$\|\mathcal{D}^s(f_m^-)\|_{F_\alpha^p} + \|\mathcal{I}^s(f_m^-)\|_{F_\alpha^p} \le C\|f\|_{F_\beta^p}$$

*for $f \in F_\beta^p$.*

*Proof.* Let $0 < p \le \infty$ and $\alpha, \beta, s$ be real numbers. Let $m$ be a positive integer. First, we note that there is a constant $C_1 = C_1(m) > 0$ such that

$$\sum_{|\gamma| \le m} |\partial^\gamma f(0)| \le C_1 \sup_{|z| \le 1} |f(z)| \qquad (4.4)$$

`ball`

for all $f \in H(\mathbb{C}^n)$ by the Cauchy estimate.

Now, given $f \in F_\beta^p$, we see from the definition of $\mathcal{D}^s$ and $\mathcal{I}^s$ that

$$|\mathcal{D}^s \partial^\gamma f(0)| + |\mathcal{I}^s \partial^\gamma f(0)| \le C|\partial^\gamma f(0)|$$

for some constant $C(s, |\gamma|) > 0$ and thus

$$|\mathcal{D}^s \partial^\gamma f(0)| + |\mathcal{I}^s \partial^\gamma f(0)| \lesssim \|f\|_{L_\beta^p}$$

by (4.4) and Lemma 3.11. Accordingly, we have

$$\|\mathcal{D}^s(f_m^-)\|_{F_\alpha^p} + \|\mathcal{I}^s(f_m^-)\|_{F_\alpha^p} \lesssim \sum_{|\gamma| \le m} (|\mathcal{D}^s \partial^\gamma f(0)| + |\mathcal{I}^s \partial^\gamma f(0)|) \|z^\gamma\|_{F_\alpha^p}$$
$$\lesssim \|f\|_{F_\beta^p},$$

as asserted. The proof is complete. $\qquad \square$

Two types of weighted Fock-Sobolev spaces with the same parameters turn out to be exactly the same, which is not too surprising in view of their definitions. More interesting is the fact that they can be identified with suitable weighted Fock spaces, as in the next theorem.

`equivnorm`

**Theorem 4.2.** *Let $s$ and $\alpha$ be real numbers. Then*

$$F_{\alpha,s}^p = \widetilde{F}_{\alpha,s}^p = \begin{cases} F_{\alpha-sp}^p & \text{if } 0 < p < \infty \\ F_{\alpha-s}^\infty & \text{if } p = \infty \end{cases}$$

*with equivalent norms.*

*Proof.* We need to prove that there is a constant $C = C(p, \alpha, s) > 0$ such that

$$C^{-1}\|f\|_{F_{\alpha-sp}^p} \le \|f\|_X \le C\|f\|_{F_{\alpha-sp}^p}, \qquad f \in H(\mathbb{C}^n)$$

for both $X = F_{\alpha,s}^p$ and $X = \widetilde{F}_{\alpha,s}^p$. Note that the second inequality of the above follows from Lemma 4.1 and Theorem 3.13.

We provide a proof of the first inequality for $X = F_{\alpha,s}^p$; the proof for $X = \widetilde{F}_{\alpha,s}^p$ is similar. Let $f \in H(\mathbb{C}^n)$. First, assume $s > 0$. Thus, using the relation $\mathcal{I}^s \mathcal{D}^s f = f$, we have by Theorem 3.13 and (4.3)

$$\|f\|_{F_{\alpha-sp}^p} = \|\mathcal{I}^s \mathcal{D}^s f\|_{F_{\alpha-sp}^p} \lesssim \|\mathcal{D}^s f\|_{F_{\alpha+sp}^p} = \|f\|_{F_{\alpha,s}^p},$$

as required. Now, assume $s < 0$. Thus, using the relation $\mathcal{I}^s \mathcal{D}^s f = f_s^+$, we have again by Theorem 3.13 and (4.3)

$$\|f_s^+\|_{F_{\alpha-sp}^p} = \|\mathcal{I}^s \mathcal{D}^s f\|_{L_{\alpha-sp}^p} \lesssim \|\mathcal{D}^s f\|_{F_{\alpha+sp}^p} \leq \|f\|_{F_{\alpha,s}^p}.$$

Since $s < 0$, we also have $\|f_s^-\|_{F_{\alpha-sp}^p} \leq \|f_s^-\|_{F_\alpha^p} \leq \|f\|_{F_{\alpha,s}^p}$. This completes the proof. $\qquad\square$

We now mention a couple of consequences of Theorem 4.2. First, we have the following parameter relation to induce the same weighted Fock-Sobolev space.

 **Corollary 4.3.** *Let* $0 < p < \infty$ *and* $s_j, \alpha_j$ *be real numbers for* $j = 1, 2$. *Then the following statements hold:*

  (a) $F_{\alpha_1,s_1}^p = F_{\alpha_2,s_2}^p$ *if and only if* $\alpha_1 - \alpha_2 = p(s_1 - s_2)$;
  (b) $F_{\alpha_1,s_1}^\infty = F_{\alpha_2,s_2}^\infty$ *if and only if* $\alpha_1 - \alpha_2 = s_1 - s_2$.

Next, we observe that the most natural definition of the weighted Fock-Sobolev spaces of positive integer order in terms of full derivatives is actually the same as the one given by fractional derivatives.

 **Corollary 4.4.** *Given* $0 < p \leq \infty$, *a positive integer* $m$ *and* $\alpha$ *real, there is a constant* $C = C(p, m, \alpha) > 0$ *such that*

$$C^{-1} \|f\|_{F_{\alpha,m}^p} \leq \sum_{|\gamma| \leq m} \|\partial^\gamma f\|_{L_\alpha^p} \leq C \|f\|_{F_{\alpha,m}^p}$$

*for* $f \in H(\mathbb{C}^n)$.

*Proof.* Let $m$ be a positive integer. Fix $f \in H(\mathbb{C}^n)$. Note that the several-variable version of (3.24) with $z \frac{\partial}{\partial z}$ replaced by $\sum_{j=0}^n z_j \partial_j$ remains valid if coefficients are appropriately adjusted. Thus we have

$$|\mathcal{R}^m f(z)| = \frac{|\mathcal{D}^m f(z)|}{(1 + |z|)^m} \lesssim \sum_{|\gamma| \leq m} |\partial^\gamma f(z)|$$

for $z \in \mathbb{C}^n$. This yields the first inequality of the corollary.

For the second inequality, we note that, given a multi-index $\gamma$, there is a constant $C_\gamma > 0$ such that

$$|\partial^\gamma K_w(z)| \leq C_\gamma (1 + |w|)^{|\gamma|} e^{\mathrm{Re}\,(z \cdot \overline{w})}, \qquad z, w \in \mathbb{C}^n;$$

recall that $K_z(w)$ denotes the Fock kernel. Thus, when $|z| \geq 1$, the estimate in Proposition 3.2 holds with $\mathcal{D}^m$ replaced by $\partial^\gamma$ for all $\gamma$ with $|\gamma| \leq m$. So, following the argument in the proof of Theorem 3.13, one obtains

$$\sum_{|\gamma| \leq m} \|\partial^\gamma f\|_{L_\alpha^p} \lesssim \|f\|_{L_{\alpha-2mp}^p}.$$

Since $\|f\|_{L^p_{\alpha-2mp}} \approx \|f\|_{F^p_{\alpha,m}}$ by Theorem 4.2, this completes the proof.          □

We now proceed to the reproducing kernels for the weighted Fock-Sobolev spaces. With Theorem 4.2 granted, we may focus on the weighted Fock spaces. The inner product on $F^2_\alpha$, inherited from $L^2_\alpha$, is given by

$$(f, g) \mapsto \int_{\mathbb{C}^n} f(z)\overline{g(z)}e^{-|z|^2}\,dV_\alpha(z).$$

However, this inner product has some disadvantage in the sense that it is not easy to find reproducing kernels explicitly. We introduce below a modified inner product (4.6), still inducing equivalent norms, which enables us to represent the weighted Fock space kernel explicitly.

It turns out that the measure

$$dW_\alpha(z) := \frac{dV(z)}{|z|^\alpha}$$

is an appropriate replacement of $dV_\alpha(z)$ to find an explicit formula for the kernel. A trouble in this case is that $|z|^{-\alpha}$ is not locally integrable near the origin when $\alpha \geq 2n$ and hence some adjustment is required. To do so we introduce the notation

$$(\psi, \varphi)_\alpha := \int_{\mathbb{C}^n} \psi(z)\overline{\varphi(z)}e^{-|z|^2}dW_\alpha(z) \qquad (4.5) \quad \boxed{\texttt{pairing}}$$

for any $\alpha$ real, whenever the integral is well defined.

Now we define an inner product $\langle\ ,\ \rangle_\alpha$ on $F^2_\alpha$ by

$$\langle f, g\rangle_\alpha := \begin{cases} (f, g)_\alpha & \text{if}\quad \alpha < 2n \\ (f^-_{\alpha/2}, g^-_{\alpha/2})_0 + (f^+_{\alpha/2}, g^+_{\alpha/2})_\alpha & \text{if}\quad \alpha \geq 2n \end{cases} \qquad (4.6) \quad \boxed{\texttt{ip}}$$

for $f, g \in F^2_\alpha$. We note from orthogonality of holomorphic monomials that, when $\alpha \geq 2n$,

$$\langle f, g\rangle_\alpha = (f, g)_\alpha \qquad (4.7) \quad \boxed{\texttt{reduce}}$$

for functions $f$ with vanishing derivatives up to order $\alpha$ at the origin. It is not hard to check that (4.6) induces an equivalent norm on $F^2_\alpha$ in case $\alpha < 2n$. Also, one may check by (4.4) and Lemma 4.1 that (4.6) induces an equivalent norm on $F^2_\alpha$ in case $\alpha \geq 2n$. So, for the rest of the paper, we will consider $F^2_\alpha$ as a Hilbert space endowed with the inner product $\langle\ ,\ \rangle_\alpha$. Also, we write $\|f\|_\alpha = \sqrt{\langle f, f\rangle_\alpha}$ for $f \in F^2_\alpha$.

Note by Proposition 2.2 that each point evaluation is a bounded linear functional on $F^2_\alpha$. So, to each $z \in \mathbb{C}^n$ there corresponds the reproducing kernel $K^\alpha_z$ such that

$$f(z) = \langle f, K^\alpha_z\rangle_\alpha$$

for $f \in F^2_\alpha$. By Proposition 2.3 holomorphic monomials span a dense subset of $F^2_\alpha$. Also, note that holomorphic monomials are mutually orthogonal in $F^2_\alpha$.

Accordingly, the set $\{z^\gamma/\|z^\gamma\|_\alpha\}_\gamma$ of normalized monomials form an orthonormal basis for $F_\alpha^2$. So, using the well-known formula

$$K^\alpha(z, w) := K_w^\alpha(z) = \sum_\gamma \phi_\gamma(z)\overline{\phi_\gamma(w)}$$

where $\{\phi_\gamma\}$ is any orthonormal basis for $F_\alpha^2$, we have

$$K^\alpha(z, w) = \sum_\gamma \frac{z^\gamma \overline{w}^\gamma}{\|z^\gamma\|_\alpha^2}. \tag{4.8}$$ `onbker`

By means of this formula, it turns out that the major part of the reproducing kernels are fractional integrals of the Fock kernel, as in the next theorem. For a more explicit formula when $\alpha$ is an even negative integer, see [4] or [5].

`ker` **Theorem 4.5.** *Let $\alpha$ be a real number. Then*

$$K^\alpha(z, w) = \mathcal{I}^{-\alpha/2} K_w(z) + E^\alpha(z, w)$$

*where the error term $E^\alpha(z, w)$ is the polynomial in $z \cdot \overline{w}$ given by*

$$E^\alpha(z, w) = \begin{cases} 0 & \text{if} \quad \alpha \leq 0 \\ \displaystyle\sum_{k \leq \alpha/2} \frac{\Gamma(n + k)}{\Gamma(n + k - \alpha/2)} \frac{(z \cdot \overline{w})^k}{k!} & \text{if} \quad 0 < \alpha < 2n \\ (K_w)_{\alpha/2}^-(z) & \text{if} \quad \alpha \geq 2n \end{cases}$$

*for $z, w \in \mathbb{C}^n$.*

*Proof.* We consider the cases $\alpha < 2n$ and $\alpha \geq 2n$ separately. First, consider the case $\alpha < 2n$. In this case an elementary computation yields

$$\|z^\gamma\|_\alpha^2 = (z^\gamma, z^\gamma)_\alpha = \frac{\gamma!\Gamma(n + |\gamma| - \alpha/2)}{\Gamma(n + |\gamma|)} \tag{4.9}$$ `mono`

for each multi-index $\gamma$. Thus, given $z, w \in \mathbb{C}^n$, a little manipulation with (4.8) yields

$$\begin{aligned} K^\alpha(z, w) &= \sum_\gamma \frac{\Gamma(n + |\gamma|)}{\Gamma(n + |\gamma| - \alpha/2)} \frac{z^\gamma \overline{w}^\gamma}{\gamma!} \\ &= \sum_{k=0}^\infty \frac{\Gamma(n + k)}{\Gamma(n + k - \alpha/2)} \frac{(z \cdot \overline{w})^k}{k!}, \end{aligned} \tag{4.10}$$ `kerseries`

which is the formula for $\alpha \leq 0$. For $0 < \alpha < 2n$ one may decompose the above sum into $\sum_{k > \alpha/2} + \sum_{k \leq \alpha/2}$ to verify the formula.

Next, consider the case $\alpha \geq 2n$. In this case (4.9) is still valid for $|\gamma| > \alpha/2$. Meanwhile, note

$$\sum_{|\gamma| \leq \alpha/2} \frac{z^\gamma \overline{w}^\gamma}{(z^\gamma, z^\gamma)_0} = \sum_{|\gamma| \leq \alpha/2} \frac{z^\gamma \overline{w}^\gamma}{\gamma!} = (K_w)_{\alpha/2}^-(z).$$

Consequently, given $z, w \in \mathbb{C}^n$, we have

$$K_w^\alpha(z) = (K_w)_{\alpha/2}^-(z) + \sum_{|\gamma| > \alpha/2} \frac{\Gamma(n + |\gamma|)}{\Gamma(n + |\gamma| - \alpha/2)} \frac{z^\gamma \overline{w}^\gamma}{\gamma!}$$

$$= (K_w)_{\alpha/2}^-(z) + \mathcal{I}^{-\alpha/2} K_w(z),$$

as asserted. This completes the proof. $\qquad\qquad\square$

By Theorem 4.5 and Proposition 3.5, we have the following estimate for the reproducing kernels.

**Corollary 4.6.** *Given* $0 < \epsilon < 1$ *and* $\alpha$ *real, there are positive constants* $C = C(\alpha, \epsilon) > 0$ *and* $\delta = \delta(\epsilon) > 0$ *such that*

$$|K^\alpha(w, z)| \leq C \times \begin{cases} 1 + |z \cdot \overline{w}|^{\alpha/2} \Lambda_{\epsilon, \delta}(z, w) & \text{if } \alpha > 0 \\ (1 + |z||w|)^{\alpha/2} \Lambda_{\epsilon, \delta}(z, w) & \text{if } \alpha \leq 0 \end{cases}$$

*for* $z, w \in \mathbb{C}^n$.

Note $|K^\alpha(z, w)| \lesssim (1 + |z||w|)^{\alpha/2} e^{|z||w|}$ by Corollary 4.6. Thus, an application of the Cauchy estimates on the ball with center $z$ and radius $1/|w|$ yields the following consequence.

**Corollary 4.7.** *Given* $\alpha$ *real and a multi-index* $\gamma$, *there is constant* $C = C(\alpha, \gamma) > 0$ *such that*

$$|\partial_z^\gamma K^\alpha(z, w)| \leq C|w|^{|\gamma|}(1 + |z||w|)^{\alpha/2} e^{|z||w|}$$

*for all* $z, w \in \mathbb{C}^n$.

Theorem 4.5 yields another consequence concerning the growth rate of the norms of the reproducing kernels. In fact, applying Corollary 4.6 and Proposition 3.8, one may verify that, given $0 < p, a < \infty$ and $\alpha$, $\beta$ real, there is a constant $C = C(p, a, \alpha, \beta) > 0$ such that

$$\int_{\mathbb{C}^n} |K^\beta(z, w)|^p e^{-a|w|^2} \, dV_\alpha(w) \leq C \frac{e^{\frac{p^2}{4a}|z|^2}}{(1 + |z|)^{\alpha - \beta p}} \tag{4.11}$$

for $z \in \mathbb{C}^n$. This immediately yields the first part of the next proposition. Recall that $F_\alpha^{\infty, 0}$ denotes the closed subspace of $F_\alpha^\infty$ defined by the condition (2.3).

**Proposition 4.8.** *Let* $0 < p \leq \infty$ *and* $\alpha, \beta$ *be real numbers. Then there is a constant* $C = C(p, \alpha, \beta) > 0$ *such that the following estimates hold for all* $w \in \mathbb{C}^n$:

(1) *For* $0 < p < \infty$,

$$\|K_w^\beta\|_{F_\alpha^p} \leq C \frac{e^{\frac{|w|^2}{2}}}{(1 + |w|)^{\alpha/p - \beta}};$$

(2) *For $p = \infty$, $K_w^\beta \in F_\alpha^{\infty,0}$ with*

$$\|K_w^\beta\|_{F_\alpha^\infty} \le C \frac{e^{\frac{|w|^2}{2}}}{(1+|w|)^{\alpha-\beta}}.$$

*Proof.* We only need to prove (2). We have $K_w^\beta \in F_\alpha^{\infty,0}$ for all $w \in \mathbb{C}^n$ by Corollary 4.7. For the norm estimate, setting

$$I(z,w) := \frac{|K_w^\beta(z)|e^{-\frac{1}{2}(|z|^2+|w|^2)}}{(1+|z|)^\alpha(1+|w|)^{\beta-\alpha}},$$

we need to show

$$\sup_{z,w\in\mathbb{C}^n} I(z,w) < \infty. \tag{4.12}$$

Using the elementary inequality

$$e^{\mathrm{Re}\,(z\cdot\overline{w})} + e^{\frac{1}{2}|z||w|} \le 2e^{\frac{1}{2}(|z|^2+|w|^2)-\frac{1}{8}|z-w|^2},$$

we have by Corollary 4.6 (with $\epsilon = 1/2$)

$$|K_w^\beta(z)| \lesssim (1+|z||w|)^{\beta/2}e^{\frac{1}{2}(|z|^2+|w|^2)-\frac{1}{8}|z-w|^2}$$

so that

$$I(z,w) \lesssim \frac{(1+|z||w|)^{\beta/2}}{(1+|z|)^\alpha(1+|w|)^{\beta-\alpha}}e^{-\frac{1}{8}|z-w|^2} \tag{4.13}$$

for all $z,w \in \mathbb{C}^n$.

To estimate the right hand side of (4.13), we consider two cases $\beta > 0$ and $\beta \le 0$ separately. When $\beta > 0$, using the inequality $(1+|z||w|) \le (1+|z|)(1+|w|)$, we have by (4.13) and (2.1)

$$I(z,w) \lesssim \left(\frac{1+|w|}{1+|z|}\right)^{\alpha-\beta/2} e^{-\frac{1}{8}|z-w|^2}$$

$$\le (1+|z-w|)^{|\alpha-\beta/2|}e^{-\frac{1}{8}|z-w|^2};$$

the constants suppressed here are independent of $z, w$. This yields (4.12) for $\beta > 0$.

Now, let $\beta \le 0$. When $|z| \ge |w|$, we have $1+|z||w| \ge 1+|w|^2 \ge (1+|w|)^2/2$ and thus

$$I(z,w) \lesssim \left(\frac{1+|w|}{1+|z|}\right)^\alpha e^{-\frac{1}{8}|z-w|^2}.$$

Similarly, when $|z| \le |w|$, we have

$$I(z,w) \lesssim \left(\frac{1+|w|}{1+|z|}\right)^{\alpha-\beta} e^{-\frac{1}{8}|z-w|^2}.$$

So, as in the case of $\beta > 0$, we conclude (4.12) for $\beta \le 0$. The proof is complete. $\square$

We now close the section by observing that a given reproducing kernel actually reproduces functions in any weighted Fock space.

**Proposition 4.9.** *Given $\alpha$ and $\beta$ real, the reproducing property*

$$f(z) = \langle f, K_z^\alpha \rangle_\alpha, \qquad z \in \mathbb{C}^n$$

*holds for $f \in F_\beta^p$ with $0 < p \leq \infty$.*

*Proof.* Fix $\alpha$ and $\beta$. Note from (2.11) that $F_\beta^p \subset F_{\beta/p}^{\infty,0}$ for any $0 < p < \infty$. Also, note $F_\beta^\infty \subset F_{\beta'}^{\infty,0}$ for $\beta' > \beta$. Thus it suffices to show that the proposition for the space $F_\beta^{\infty,0}$. Given $z \in \mathbb{C}^n$, we claim that there is a constant $C_z = C_z(\alpha, \beta) > 0$ such that

$$|\langle f, K_z^\alpha \rangle_\alpha| \leq C_z \|f\|_{F_\beta^\infty} \tag{4.14}$$

for $f \in F_\beta^\infty$. With this granted, we conclude the asserted reproducing property for the space $F_\beta^{\infty,0}$, because holomorphic polynomials form a dense subset in that space by Proposition 2.3.

It remains to show (4.14). Let $f \in F_\beta^\infty$. The case $\alpha \leq 0$ is easily handled, because

$$|\langle f, K_z^\alpha \rangle_\alpha| = |(f, K_z^\alpha)_\alpha| \leq \|f\|_{F_\beta^\infty} \|K_z^\alpha\|_{F_{\alpha-\beta}^1}.$$

Now, assume $\alpha > 0$. Since $K_z^\alpha$ reproduces holomorphic polynomials, we have

$$\langle f, K_z^\alpha \rangle_\alpha = \langle f_\alpha^+, K_z^\alpha \rangle_\alpha + \langle f_\alpha^-, K_z^\alpha \rangle_\alpha = (f_\alpha^+, K_z^\alpha)_\alpha + f_\alpha^-(z)$$

by (4.7) even when $\alpha \geq 2n$.

Note $\|f_\alpha^-\|_{F_\beta^\infty} \lesssim \|f\|_{F_\beta^\infty}$ by (2.10). Thus we have $|f_\alpha^-(z)| \leq C_z \|f\|_{F_\beta^\infty}$ by Proposition 2.2 and $\|f_\alpha^+\|_{F_\beta^\infty} \lesssim \|f\|_{F_\beta^\infty}$. Also, note

$$|(f_\alpha^+, K_z^\alpha)_\alpha| \leq \|g\|_{L_\beta^\infty} \|K_z^\alpha\|_{F_{-\beta}^1}$$

where $g(w) = |f_\alpha^+(w)| |w|^{-\alpha}$. Note $|g(w)| \lesssim \|f\|_{F_\beta^\infty}$ for $|w| \leq 1$ by (2.12). Accordingly, $\|g\|_{L_\beta^\infty} \lesssim \|f\|_{F_\beta^\infty} + \|f_\alpha^+\|_{F_\beta^\infty} \lesssim \|f\|_{F_\beta^\infty}$. So, we obtain (4.14). The proof is complete. □

## 5. Applications

In this section we apply the results obtained in earlier sections to derive some basic properties of the Fock-Sobolev spaces such as projections, dual spaces, complex interpolation spaces and Carleson measures. Those were first studied by Cho and Zhu [5] when the Sobolev order is a positive integer. Here, we extend their results to an arbitrary order. In fact our results, even when restricted to an order of positive integer, contain their results as special cases (except for Carleson measures). For the extension to an arbitrary order, note that Theorem 4.2 allows us to focus on the weighted Fock spaces throughout the section.

In addition to the results we have established so far, we need some additional technical preliminaries. We begin with by recalling the reproducing property

$$f(z) = (f, K_z^\alpha)_\alpha \quad \text{for} \quad \alpha < 2n$$

and for any weighted Fock-function $f$. Also, introducing the truncated kernel

$$K_w^{\alpha,+}(z) = K^{\alpha,+}(z,w) := (K_w^\alpha)_\alpha^+(z),$$

we have by (4.7) the reproducing property

$$f_\alpha^+(z) = (f_\alpha^+, K_z^\alpha)_\alpha = (f, K_z^{\alpha,+})_\alpha \quad \text{for} \quad \alpha \geq 2n \tag{5.1} \quad \boxed{\texttt{rep+}}$$

and for any weighted Fock-function $f$. Motivated by these reproducing kernels, we first consider auxiliary integral operators $S_\alpha$ and $S_\alpha^+$ defined by

$$S_\alpha \psi(z) := \int_{\mathbb{B}_n} \psi(w)|K^\alpha(z,w)|e^{-|w|^2}\, dW_\alpha(w) \quad \text{for} \quad \alpha < 2n$$

and

$$S_\alpha^+ \psi(z) := \int_{\mathbb{B}_n} \psi(w)|K^{\alpha,+}(z,w)|e^{-|w|^2}\, dW_\alpha(w) \quad \text{for} \quad \alpha \geq 2n;$$

recall that $\mathbb{B}_n$ denotes the unit ball of $\mathbb{C}^n$.

$\boxed{\texttt{sabdd}}$ **Lemma 5.1.** *Given $\beta$ real and $1 \leq p \leq \infty$ the following statements hold:*

(1) *If $\alpha < 2n$, then $S_\alpha : L_\beta^p \to L_\beta^p$ is bounded;*
(2) *If $\alpha \geq 2n$, then $S_\alpha^+ : L_\beta^p \to L_\beta^p$ is bounded.*

*Proof.* We provide the details for $\alpha \geq 2n$. In case $\alpha < 2n$, one may easily modify the proof below, because $|w|^{-\alpha}$ is integrable near the origin.

Fix any real number $\beta$ and let $\alpha \geq 2n$. Given any $\beta'$ real, we have by (2.12) and Proposition 4.8

$$\sup_{w \in \mathbb{B}_n} \frac{|K^{\alpha,+}(z,w)|}{|w|^\alpha} \leq C\|K_z^\alpha\|_{F_{\beta'}^\infty} \leq C\frac{e^{\frac{1}{2}|z|^2}}{(1+|z|)^{\beta'-\alpha}}, \qquad z \in \mathbb{C}^n \tag{5.2} \quad \boxed{\texttt{beta'}}$$

for some constant $C = C(\alpha, \beta') > 0$. Thus, choosing $\beta' = \alpha - \beta + 2n + 1$, we have

$$|S_\alpha^+ \psi(z)| \lesssim \frac{e^{\frac{1}{2}|z|^2}}{(1+|z|)^{2n+1-\beta}}\|\psi\|_{L_\beta^1}.$$

This implies that $S_\alpha^+$ is bounded on $L_\beta^1$. Also, choosing $\beta' = \alpha - \beta$, we obtain

$$|S_\alpha^+ \psi(z)| \lesssim \frac{e^{\frac{1}{2}|z|^2}}{(1+|z|)^{-\beta}}\|\psi\|_{L_\beta^\infty}.$$

So, $S_\alpha^+$ is bounded on $L_\beta^\infty$. In particular, $S_\alpha^+$ is bounded on $L_0^\infty$. Thus, $S_\alpha^+$ is also bounded on $L_\beta^p$ for any $1 < p < \infty$ by the Stein interpolation theorem. The proof is complete. $\qquad \square$

Next, we introduce a class of auxiliary function spaces and a related integral operator. For $r > 0$, let $\Omega_r := \mathbb{C}^n \setminus r\mathbb{B}_n$. For $0 < p < \infty$ and $\alpha$ real, we denote by $\mathcal{L}_\alpha^{p,r} = \mathcal{L}_\alpha^{p,r}(\Omega_r)$ the space of all Lebesgue measurable functions $\psi$ on $\Omega_r$ such that the norm

$$\|\psi\|_{\mathcal{L}_\alpha^{p,r}} := \left\{ \int_{\Omega_r} \left| \psi(z) e^{-\frac{1}{2}|z|^2} \right|^p dW_\alpha(z) \right\}^{1/p}$$

is finite. When $p = \infty$, we denote by $\mathcal{L}_\alpha^{\infty,r} = \mathcal{L}_\alpha^{\infty,r}(\Omega_r)$ the space of all Lebesgue measurable functions $\psi$ on $\Omega_r$ such that the norm

$$\|\psi\|_{\mathcal{L}_\alpha^{\infty,r}} := \operatorname{esssup} \left\{ \frac{|\psi(z)| e^{-\frac{1}{2}|z|^2}}{|z|^\alpha} : z \in \Omega_r \right\}$$

is finite. Note

$$\mathcal{L}_\alpha^{p,r} \cap H(\mathbb{C}^n) = F_\alpha^p \tag{5.3}$$ `identify`

with equivalent norms. Here, we are identifying an entire function with its restriction to $\Omega_r$.

For $\alpha$ real and $r > 0$, consider an integral operator $T_\alpha$ defined by

$$T_\alpha^r \psi(z) := \int_{\Omega_r} \psi(w)(1 + |z||w|)^{\alpha/2} \Lambda(z,w) e^{-|w|^2} dW_\alpha(w), \qquad z \in \mathbb{C}^n$$

where $\Lambda := \Lambda_{\epsilon,\delta}$ denotes any (fixed) function as in Corollary 4.6.

`stbdd` **Lemma 5.2.** *Given $r > 0$ and $\alpha$ real, the operator $T_\alpha^r : \mathcal{L}_\beta^{p,r} \to L_\beta^p$ is bounded for any $\beta$ real and $1 \le p \le \infty$.*

*Proof.* Fix $r > 0$ and real numbers $\alpha, \beta$. First, assume $\alpha > 0$. Since $(1+|z||w|) \le (1+|z|)(1+|w|)$, we have

$$|T_\alpha^r \psi(z)| \le \int_{\Omega_r} |\psi(w)|(1+|z|)^{\alpha/2}(1+|w|)^{\alpha/2} \Lambda(z,w) e^{-|w|^2} dW_\alpha(w)$$

$$\approx \int_{\Omega_r} |\psi(w)| \left( \frac{1+|z|}{1+|w|} \right)^{\alpha/2} \Lambda(z,w) e^{-|w|^2} dV(w).$$

Thus, denoting by $\widetilde{\psi}$ denotes the extension of $\psi$ defined to be $0$ on $r\mathbb{B}_n$, we obtain

$$|T_\alpha^r \psi(z)| \lesssim L_{\alpha/2} \widetilde{\psi}(z) \tag{5.4}$$ `tapz`

where $L_{\alpha/2}$ is one of the operators considered in Proposition 3.9. Now, since $\|\widetilde{\psi}\|_{L_\beta^p} \approx \|\psi\|_{\mathcal{L}_\beta^{p,r}}$, we conclude by Proposition 3.9 that $T_\alpha^r : \mathcal{L}_\beta^{p,r} \to L_\beta^p$ is bounded for any $1 \le p \le \infty$.

Next, assume $\alpha \le 0$. In this case, note $1 + |z||w| \approx |z||w| \approx (1+|z|)(1+|w|)$ for $z, w \in \Omega_r$. Thus, the argument above shows that (5.4) remains valid for $z \in \Omega_r$. On the other hand, note $\Lambda(z,w) \le 2e^{|z||w|} \le 2e^{r|w|}$ for $z \in r\mathbb{B}_n$. Also, note $e^{r|w|}$ belongs to $\mathcal{L}_s^{q,r}$ for any $s$ real and $1 \le q \le \infty$. Thus, given $1 \le p \le \infty$, we have by Hölder's inequality

$$\sup_{z \in r\mathbb{B}_n} |T_\alpha^r \psi(z)| \le C \|\psi\|_{\mathcal{L}_\beta^{p,r}}$$

for some constant $C = C(p, \beta, r) > 0$. Accordingly, we obtain

$$|T_\alpha^r \psi| \lesssim \|\psi\|_{\mathcal{L}_\beta^{p,r}} + L_{\alpha/2}\widetilde{\psi}.$$

This implies as before that $T_\alpha^r : \mathcal{L}_\alpha^{p,r} \to L_\beta^p$ is bounded. The proof is complete. $\quad\square$

`projection`

5.1. **Reproducing operator.** Assume $\alpha \geq 2n$ for a moment. Let $H_\alpha(\mathbb{C}^n)$ be the class of all entire functions $f \in H(\mathbb{C}^n)$ such that $f_\alpha^- = 0$. For $\beta$ real, put

$$F_\beta^{p,\alpha+} := F_\beta^p \cap H_\alpha(\mathbb{C}^n), \qquad 0 < p \leq \infty,$$

which is regraded as a closed subspace of $F_\beta^p$. Note from (5.1) that $K_z^{\alpha,+}$ is the reproducing kernel at $z$ for the spaces $F_\beta^{p,\alpha+}$ under the pairing $(\ ,\ )_\alpha$.

We now introduce integral operators induced by the reproducing kernels. Namely, we define

$$P_\alpha \psi(z) := (\psi, K_z^\alpha)_\alpha \quad \text{for} \quad \alpha < 2n$$

and

$$P_\alpha^+ \psi(z) := (\psi, K_z^{\alpha,+})_\alpha \quad \text{for} \quad \alpha \geq 2n.$$

Associated with these operators are the operators

$$Q_\alpha \psi(z) := (\psi, |K_z^\alpha|)_\alpha \quad \text{for} \quad \alpha < 2n$$

and

$$Q_\alpha^+ \psi(z) := (\psi, |K_z^{\alpha,+}|)_\alpha \quad \text{for} \quad \alpha \geq 2n.$$

Note $|P_\alpha \psi| \leq Q_\alpha |\psi|$ and $|P_\alpha^+ \psi| \leq Q_\alpha^+ |\psi|$.

`pabdd`

**Theorem 5.3.** *Given $\beta$ real and $1 \leq p \leq \infty$ the following statements hold:*

(1) *If $\alpha < 2n$, then $Q_\alpha : L_\beta^p \to L_\beta^p$ is bounded and $P_\alpha : L_\beta^p \to F_\beta^p$ is a bounded projection;*

(2) *If $\alpha \geq 2n$, then $Q_\alpha^+ : L_\beta^p \to L_\beta^p$ is bounded and $P_\alpha^+ : L_\beta^p \to F_\beta^{p,\alpha+}$ is a bounded projection.*

*Proof.* We provide the details for (2); the proof for (1) is similar. Fix $\beta$ and $1 \leq p \leq \infty$. Let $\alpha \geq 2n$. Since $|K_z^{\alpha,+}(w)| \leq |K_z^\alpha(w)| + |(K_z^\alpha)_\alpha^-(w)|$ and $|(K_z^\alpha)_\alpha^-(w)| \lesssim 1 + |z \cdot \overline{w}|^\alpha$, it is easily seen that $K_z^{\alpha,+}(w)$ satisfies the same growth estimate given for the original kernel in Corollaries 4.6. Thus we have

$$|Q_\alpha^+ \psi| \lesssim S_\alpha^+ |\psi| + T_\alpha^1 |\widetilde{\psi}|$$

where $\widetilde{\psi}$ is the restriction of $\psi$ to $\Omega_1$. Since $\|\widetilde{\psi}\|_{\mathcal{L}_\beta^{p,1}} \lesssim \|\psi\|_{L_\beta^p}$, the above, together with Lemmas 5.1 and 5.2, implies that $Q_\alpha^+ : L_\beta^p \to L_\beta^p$ is bounded.

Next, $P_\alpha^+$ clearly has the same boundedness properties as $Q_\alpha^+$, because $|P_\alpha^+ \psi| \leq Q_\alpha^+ |\psi|$. Note that the estimate in Corollary 4.7 with $K^\alpha$ replaced by $K^{\alpha,+}$ remains valid. Thus, one can justify the differentiation under the integral sign to verify that $P_\alpha^+$ takes $L_\beta^p$ into $H_\alpha(\mathbb{C}^n)$. So, $P_\alpha^+ : L_\beta^p \to F_\beta^{p,\alpha+}$ is a bounded projection by the reproducing property. The proof is complete. $\quad\square$

*Remark.* Let $s \geq 0$ and $\alpha$ be a real number. Further assume $s \geq \alpha$ if $\alpha \geq 2n$. For such $s$ and $\alpha$, consider an integral operator

$$P_{\alpha,s}\psi(z) := \left(\psi, (K_z^\alpha)_s^+\right)_\alpha.$$

For example, $P_{\alpha,\alpha} = P_\alpha^+$ when $\alpha \geq 2n$. One may easily deduce from Theorem 5.3 that $P_{\alpha,s} : L_\beta^p \to F_\beta^{p,s+}$ is a bounded projection for any $\beta$ real.

In connection with the reproducing operators considered above, one may naturally consider the integral operators on the spaces $\mathcal{L}_\alpha^{p,r}$ induced by the reproducing kernel. More explicitly, given $\alpha$ real and $r > 0$, the integral pairing

$$[\psi, \varphi]_\alpha^r := \int_{\Omega_r} \psi(w)\overline{\varphi(w)}e^{-|w|^2}\,dW_\alpha(w)$$

gives rise to an integral operator

$$P_\alpha^r\psi(z) := [\psi, K_z^\alpha]_\alpha^r.$$

When $\alpha < 2n$, note that $P_\alpha$ is formally the limiting operator of $P_\alpha^r$ as $r \to 0^+$. So, in view of Theorem 5.3, one may expect, especially when $r$ is small, some useful property for those operators. Such an intuition is formulated in the next theorem, which is one of the keys to our results on duality and complex interpolation later.

**Theorem 5.4.** *Let $1 \leq p \leq \infty$, $r > 0$ and $\alpha$, $\beta$ be real numbers. Then the operator $P_\alpha^r : \mathcal{L}_\beta^{p,r} \to F_\beta^p$ is bounded. Moreover, $P_\alpha^r : F_\beta^p \to F_\beta^p$ is invertible for all $r$ sufficiently small (depending on $\alpha$, $\beta$ and $p$).*

*Proof.* Note $|P_\alpha^r\psi| \lesssim T_\alpha^r|\psi|$ by Corollary 4.6. Thus $P_\alpha^r : \mathcal{L}_\beta^{p,r} \to \mathcal{L}_\beta^{p,r}$ is bounded by Lemma 5.2. Also, $P_\alpha^r$ takes $\mathcal{L}_\beta^{p,r}$ into $H(\mathbb{C}^n)$, as in the proof of Theorem 5.3. Thus the first part holds by (5.3).

For the second part, we provide details only for the case $\alpha \geq 0$; the case $\alpha < 0$ is similar and simpler. Let $X_\alpha$ be the set of all holomorphic polynomials of degree at most $\alpha$. Also, let $Y_{\beta,\alpha}^p$ be the closed subspace of $F_\beta^p$ consisting of all functions with vanishing derivatives at the origin up to order $\alpha$. Note that $P_\alpha^r$ maps a monomial to another monomial of the same multi-degree. Accordingly, $P_\alpha^r : X_\alpha \to X_\alpha$ is invertible. Also, $P_\alpha^r$ takes $Y_{\beta,\alpha}^p$ boundedly into itself. So, it suffices to show that $P_\alpha^r : Y_{\beta,\alpha}^p \to Y_{\beta,\alpha}^p$ is invertible for all $r$ sufficiently small.

To begin with, let $r \leq 1$. Let $f \in Y_{\beta,\alpha}^p$. Since $f(z) = (f, K_z^\alpha)_\alpha$ by (5.1), we have

$$f(z) - P_\alpha^r f(z) = \int_{r\mathbb{B}_n} f(w)K^\alpha(z,w)e^{-|w|^2}\,dW_\alpha(w).$$

Now, it follows from Corollary 4.7 and (2.12) that

$$|f(z) - P_\alpha^r f(z)| \lesssim r^{2n}(1+|z|)^{\alpha/2}e^{|z|}\left(\sup_{w \in r\mathbb{B}_n} \frac{|f(w)|}{|w|^\alpha}\right)$$

$$\lesssim r^{2n}(1+|z|)^{\alpha/2}e^{|z|}\|f\|_{F_\beta^p}$$

inverse

for all $z \in \mathbb{C}^n$. This yields

$$\|f - P_\alpha^r f\|_{F_\beta^p} \lesssim r^{2n} \|f\|_{F_\beta^p};$$

the constant suppressed here is independent of $r$ and $f$. It follows that

$$\|I - P_\alpha^r\| \lesssim r^{2n} \to 0 \quad \text{as} \quad r \to 0^+$$

where $I$ is the identity operator on $Y_{\beta,\alpha}^p$ and $\|I - P_\alpha^r\|$ is the operator norm of $I - P_\alpha^r$ acting on $Y_{\beta,\alpha}^p$. So, $P_\alpha^r : Y_{\beta,\alpha}^p \to Y_{\beta,\alpha}^p$ is invertible for all $r$ sufficiently small. The proof is complete. $\qquad\square$

For $\alpha \geq 2n$, we remark that the analogue of Theorem 5.4 for the operator

$$P_\alpha^{r,+} \psi(z) := [\psi, K_z^{\alpha,+}]_\alpha^r$$

is also true by a similar argument.

5.2. **Duality.** We identify the dual of the weighted Fock spaces under the suitably chosen pairing depending on the parameters $\alpha$. In what follows the superscript $*$ stands for the dual of the underlying space.

First, we consider the case $1 < p < \infty$. Before proceeding, we note the duality $(\mathcal{L}_{\alpha p}^{p,r})^* = \mathcal{L}_{\beta q}^{q,r}$ under the pairing $[\ ,\ ]_{\alpha+\beta}^r$ where $1/p + 1/q = 1$.

**Theorem 5.5.** *Let $1 < p < \infty$ and $\alpha$, $\beta$ be real numbers. Then*

$$(F_{\alpha p}^p)^* = F_{\beta q}^q$$

*with equivalent norms under the pairing $\langle\ ,\ \rangle_{\alpha+\beta}$. Here, $q$ denotes the conjugate index of $p$.*

*Proof.* Put $s := (\alpha + \beta)/2$ for short. We provide details for the case $s \geq n$; the case $s < n$ is similar and simpler. Let $f \in F_{\alpha p}^p$ and $g \in F_{\beta q}^q$. Then we have by Hölder's inequality

$$|\langle f, g \rangle_{2s}| \leq |(f_s^+, g_s^+)_{2s}| + |(f_s^-, g_s^-)_0|$$
$$\leq \|f_s^+\|_{F_{2p}^p}\|g_s^+\|_{F_{2q}^q} + \|f_s^-\|_{F_0^p}\|g_s^-\|_{F_0^q}.$$

This, together with Lemma 4.1, yields $|\langle f, g \rangle_{2s}| \lesssim \|f\|_{F_{\alpha p}^p}\|g\|_{F_{\beta q}^q}$. Consequently, we see that $F_{\beta q}^q$ is continuously embedded into $(F_{\alpha p}^p)^*$ via the given pairing.

Conversely, let $\nu$ be a bounded linear functional on $F_{\alpha p}^p$. Pick any $r > 0$, say $r = 1$. Then, according to the Hahn-Banach extension theorem, $\nu$ can be extended (without increasing its norm) to a bounded linear functional on $\mathcal{L}_{\alpha p}^{p,1}$. Using the duality $(\mathcal{L}_{\alpha p}^{p,1})^* = \mathcal{L}_{\beta q}^{q,1}$ under the pairing $[\ ,\ ]_{2s}^1$, we can pick some $\psi \in \mathcal{L}_{\beta q}^{q,1}$ such that $\|\psi\|_{\mathcal{L}_{\beta q}^{q,1}} \leq \|\nu\|$ and $\nu(f) = [f, \psi]_{2s}^1$ for all $f \in F_{\alpha p}^p$. Put $g := P_{2s}^1 \psi$. Note $g \in F_{\beta q}^q$ with $\|g\|_{F_{\beta q}^q} \lesssim \|\nu\|$ by Theorem 5.4. We claim

$$\nu(f) = \langle f, g \rangle_{2s} \qquad (5.5)$$

for all $f \in F_{\alpha p}^p$. With this granted, we conclude that $(F_{\alpha p}^p)^*$ is continuously embedded into $F_{\beta q}^q$, as required.

To prove (5.5), fix any $f \in F_{\alpha p}^p$. Since $(K_z^{2s})_s^+ = \mathcal{I}^{-s} K_z$ and $(K_z^{2s})_s^- = (K_z)_s^-$ by Theorem 4.5, we have by the reproducing property

$$f(z) = \langle f, K_z^{2s} \rangle_{2s} = \left( f_s^+, \mathcal{I}^{-s} K_z \right)_{2s} + \left( f_s^-, (K_z)_s^- \right)_0.$$

We also have

$$g_s^+(w) = \left[ \psi, \mathcal{I}^{-s} K_w \right]_{2s}^1 \quad \text{and} \quad g_s^-(w) = \left[ \psi, (K_w)_s^- \right]_{2s}^1.$$

So, we obtain

$$\begin{aligned}
\nu(f) &= [f, \psi]_{2s}^1 \\
&= \left[ \left( f_s^+, \mathcal{I}^{-s} K_z \right)_{2s}, \psi(z) \right]_{2s}^1 + \left[ \left( f_s^-, (K_z)_s^- \right)_0, \psi(z) \right]_{2s}^1 \\
&= \left( f_s^+(w), \left[ \psi, \mathcal{I}^{-s} K_w \right]_{2s}^1 \right)_{2s} + \left( f_s^-, \left[ \psi, (K_w)_s^- \right]_{2s}^1 \right)_0 \\
&= \langle f, g \rangle_{2s},
\end{aligned}$$

as claimed. In the third equality of the above the applications of Fubini's theorem are justified by Proposition 3.5 and the boundedness of $T_{2s}^1 : \mathcal{L}_{\beta q}^{q,1} \to L_{\beta q}^q$ (Lemma 5.2). The proof is complete. $\qquad\square$

For the duality when $0 < p \leq 1$, we need the following density property of the reproducing kernels.

**Lemma 5.6.** *Given $\alpha$ and $\beta$ real, the span of $\{K_w^\alpha : w \in \mathbb{C}^n\}$ is dense in $F_\beta^{\infty,0}$ and $F_\beta^p$ for any $0 < p < \infty$.*

*Proof.* In the special case when $\alpha = \beta$ is a negative even integer, this is proved in [5, Lemma 17]. Following the idea of the proof of [5, Lemma 17], we introduce auxiliary Hilbert function spaces $_tF_0^2 := L^2(e^{-t|z|^2}dV) \cap H(\mathbb{C}^n)$ for $t > 0$. Note from Lemma 2.1 that there is a constant $C = C(t) > 0$ such that

$$|f(z)| \leq Ce^{t|z|^2/2}\|f\|_{_tF_0^2}, \qquad z \in \mathbb{C}^n \tag{5.6}$$

for all $f \in {_tF_0^2}$.

Fix any $\alpha$ real. Note that the kernel functions $K_w^\alpha$ are all contained in $_tF_0^2$ by Corollary 4.7. We first show that the span of $\{K_w^\alpha : w \in \mathbb{C}^n\}$ is dense in $_tF_0^2$. Let $f \in {_tF_0^2}$ and assume

$$\int_{\mathbb{C}^n} f(z) K^\alpha(w, z) e^{-t|z|^2} \, dV(z) = 0$$

for all $w \in \mathbb{C}^n$. Note from Theorem 4.5 that $K^\alpha(w, z)$ is a series in $w \cdot \overline{z}$ with nonzero coefficients. So, differentiating at the origin as many times as needed under the integral sign in the left hand side of the above, which is justified by Corollary 4.7, we see that $f$ is orthogonal to all holomorphic monomials. Since holomorphic monomials span a dense subset of $_tF_0^2$, we conclude $f = 0$ and thus that the span of $\{K_w^\alpha : w \in \mathbb{C}^n\}$ is dense in $_tF_0^2$.

Fix any $\beta$ real and $0 < p \leq \infty$. Also, fix $0 < t < 1$. Let $g \in H(\mathbb{C}^n)$ be an arbitrary polynomial. For any function $h$ in the span of $\{K_w^\alpha : w \in \mathbb{C}^n\}$, we have

by (5.6)

$$\|g - h\|_{F_\beta^p} \lesssim \|g - h\|_{tF_0^2};$$

the constant suppressed here is independent of $g$ and $h$. Since the span of $\{K_w^\alpha : w \in \mathbb{C}^n\}$ is dense in $_tF_0^2$, one can make the right hand side of the above arbitrarily small by choosing suitable $h$. Now, since the set of all holomorphic polynomials is dense in the space under consideration by Proposition 2.3, we conclude the lemma. $\quad\square$

**Theorem 5.7.** *Let* $0 < p \leq 1$ *and* $\alpha, \beta$ *be real numbers. Then*

$$(F_{\alpha p}^p)^* = F_\beta^\infty$$

*with equivalent norms under the pairing* $\langle\ ,\ \rangle_{\alpha+\beta}$.

*Proof.* Put $s := (\alpha + \beta)/2$ again for short. As in the proof of Theorem 5.5, we provide details only for the case $s \geq n$. Let $f \in F_{\alpha p}^p$ and $g \in F_\beta^\infty$. We have by Hölder's inequality and Proposition 3.10

$$|(f_s^-, g_s^-)_0| \lesssim \|f_s^-\|_{F_0^p} \|g_s^-\|_{F_0^\infty}.$$

On the other hand, we have

$$|(f_s^+, g_s^+)_{2s}| \leq \int_{\mathbb{C}^n} |f_s^+(z)||g_s^+|e^{-|z|^2}\, dW_{2s}(z)$$

$$\leq \sup_{|z| \leq 1} \frac{|f_s^+(z)||g_s^+(z)|}{|z|^{2s}} + \|g_s^+\|_{F_\beta^\infty} \int_{|z| \geq 1} |f_s^+(z)|e^{-\frac{1}{2}|z|^2}\, dW_\alpha(z).$$

By (2.12) the first term of the above is dominated by some constant times $\|f\|_{F_{\alpha p}^p}\|g\|_{F_\beta^\infty}$. The integral in the second term is comparable to

$$\int_{\mathbb{C}^n} |f_s^+(z)|e^{-\frac{1}{2}|z|^2}\, dV_\alpha(z) \lesssim \|f_s^+\|_{F_{\alpha p}^p};$$

the last estimate comes from Proposition 3.10. Thus we obtain

$$|\langle f, g\rangle_{2s}| \lesssim \|f\|_{F_{\alpha p}^p}\|g\|_{F_\beta^\infty} + \|f_s^+\|_{F_{\alpha p}^p}\|g_s^+\|_{F_\beta^\infty} + \|f_s^-\|_{F_0^p}\|g_s^-\|_{F_0^\infty}$$

$$\lesssim \|f\|_{F_{\alpha p}^p}\|g\|_{F_\beta^\infty};$$

the last estimate comes from Lemma 4.1. This shows that $F_\beta^\infty$ is continuously embedded in $(F_{\alpha p}^p)^*$ via the integral pairing $\langle\ ,\ \rangle_{2s}$.

Conversely, let $\nu$ be a bounded linear functional on $F_{\alpha p}^p$. Put

$$g(w) = \overline{\nu(K_w^{2s})}, \qquad w \in \mathbb{C}^n. \tag{5.7}$$

This $g$ is well defined, because $K_w^{2s} \in F_{\alpha p}^p$ with

$$\|K_w^{2s}\|_{F_{\alpha p}^p} \lesssim (1 + |w|)^\beta e^{\frac{|w|^2}{2}}. \tag{5.8}$$

by Proposition 4.8. Note from Theorem 4.5 that the homogeneous expansion of $K_w^{2s}$ is given by

$$K_w^{2s}(z) = \sum_{k \le s} \frac{(z \cdot \overline{w})^k}{k!} + \sum_{k > s} \frac{\Gamma(n+k)}{\Gamma(n+k-s)} \frac{(z \cdot \overline{w})^k}{k!}. \tag{5.9}$$

<div align="right"><code>kerhomoexp</code></div>

By a direct computation via Sterling's formula one can check

$$\frac{\big\| \|z|^k \big\|_{L^p_{\alpha p}}^p}{(k!)^{p/2}} \approx k^{n-1/2-p(1+2\alpha)/4} \quad \text{as} \quad k \to \infty.$$

Using the inequality $|z \cdot \overline{w}| \le |z||w|$ and the above estimate, one can see that the series (5.9) converges in the norm topology of $F_{\alpha p}^p$ whenever $w$ is restricted to a compact subset of $\mathbb{C}^n$. Thus the series

$$\sum_{k > s} \frac{\Gamma(n+k)}{\Gamma(n+k-s)} \frac{\overline{\nu\big[\big((\cdot) \cdot \overline{w}\big)^k\big]}}{k!}$$

converges uniformly on compacta, which clearly implies that $g$ is entire. In addition, since $|g(w)| \le \|\nu\| \|K_w^{2s}\|_{F_{\alpha p}^p}$ by the boundedness of $\nu$ on $F_{\alpha p}^p$, we have by (5.8) $g \in F_\beta^\infty$ with $\|g\|_{F_\beta^\infty} \lesssim \|\nu\|$. We also have by the reproducing property

$$\nu(K_w^{2s}) = \overline{g(w)} = \overline{\langle g, K_w^{2s} \rangle_{2s}} = \langle K_w^{2s}, g \rangle_{2s}.$$

Since this is true for all $w \in \mathbb{C}^n$, we conclude $\nu = \langle \, \cdot \, , g \rangle_{2s}$ by Lemma 5.6. This completes the proof. $\square$

<div><code>littledual</code></div>

**Theorem 5.8.** *Let $\alpha$ and $\beta$ be real numbers. Then*

$$\left( F_\beta^{\infty,0} \right)^* = F_\alpha^1$$

*with equivalent norms under the pairing $\langle \, , \, \rangle_{\alpha+\beta}$.*

*Proof.* This time we put $s := \alpha + \beta$. We see from Theorem 5.7 that $F_\alpha^1$ is continuously embedded in $\left( F_\beta^{\infty,0} \right)^*$ via the integral pairing $\langle \, , \, \rangle_s$.

Conversely, let $\nu$ be a bounded linear functional on $F_\beta^{\infty,0}$. Note $K_w^s \in F_\beta^{\infty,0}$ by Proposition 4.8. So, setting $g(w) := \overline{\nu(K_w^s)}$ and proceeding as in the proof of Theorem 5.7, we see that $g$ is entire. Note from Proposition 4.8

$$|g(w)| \le \|\nu\| \|K_w^s\|_{F_\beta^\infty} \lesssim e^{\frac{1}{2}|w|^2}(1+|w|)^\alpha \|\nu\| \tag{5.10}$$

<div align="right"><code>gw</code></div>

for all $w$.

We now proceed to show

$$g \in F_\alpha^1 \quad \text{with} \quad \|g\|_{F_\alpha^1} \le C\|\nu\| \tag{5.11}$$

<div align="right"><code>gL1</code></div>

for some constant $C > 0$ independent of $\nu$ and $g$. For $f \in F_\beta^\infty$ and $0 < r < 1$, let $f_r$ be the dilated function $z \mapsto f(rz)$. Then $f_r \in F_\beta^{\infty,0}$ and

$$\|f_r\|_{F_\beta^\infty} \le r^{-|\beta|}\|f\|_{F_\beta^\infty}$$

for all $0 < r < 1$. Thus, defining $\nu_r(f) := \nu(f_r)$, we see that $\nu_r \in (F_\beta^\infty)^*$ with $\|\nu_r\| \leq r^{-|\beta|}\|\nu\|$. Note from (5.10) that $g_r \in F_\alpha^1$. Also, note $(K_w^s)_r = K_{rw}^s$. Thus we have

$$\nu_r(K_w^s) = \nu(K_{rw}^s) = \overline{g_r(w)} = \langle K_w^s, g_r \rangle_s$$

for all $w$. It follows from the above and Lemma 5.6 that $\nu_r(h) = \langle h, g_r \rangle_s$ for $h \in F_\beta^{\infty,0}$. In particular, we have

$$\nu_r(f_r) = \langle f_r, g_r \rangle_s = \langle f, g_{r^2} \rangle_s;$$

the last equality holds by the orthogonality of monomials. It follows that

$$|\langle f, g_{r^2} \rangle_s| \leq \|\nu_r\|\|f_r\|_{F_\beta^\infty} \leq r^{-2|\beta|}\|\nu\|\|f\|_{F_\beta^\infty} \qquad (5.12) \quad \boxed{\text{fgr}}$$

for $f \in F_\beta^\infty$.

Note $P_s^1 : \mathcal{L}_\beta^{\infty,1} \to F_\beta^\infty$ is bounded by Theorem 5.4. Now, given $\psi \in C_c(\Omega_1)$, put $\widetilde{\psi}(w) = e^{\frac{1}{2}|w|^2}|w|^\beta\psi(w)$. Note $P_s^1\widetilde{\psi} \in F_\beta^\infty$ with

$$\|P_s^1\widetilde{\psi}\|_{F_\beta^\infty} \leq \|P_s^1\|\|\widetilde{\psi}\|_{\mathcal{L}_\beta^{\infty,1}} = \|P_s^1\|\|\psi\|_{L^\infty(dV)} \qquad (5.13) \quad \boxed{\text{tildepsi}}$$

where $\|P_s^1\|$ denotes the operator norm of $P_s^1 : \mathcal{L}_\beta^{\infty,1} \to F_\beta^\infty$. Meanwhile, proceeding as in the proof of Theorem 5.5, we have

$$[\widetilde{\psi}, g_{r^2}]_s^1 = \langle P_s^1\widetilde{\psi}, g_{r^2} \rangle_s$$

and thus by (5.12) and (5.13)

$$\left| \int_{|w| \geq 1} \psi(w)\overline{g_{r^2}(w)}\frac{e^{-\frac{1}{2}|w|^2}}{|w|^\alpha} \, dV(w) \right| \lesssim r^{-2|\beta|}\|\nu\|\|P_s^1\|\|\psi\|_{L^\infty(dV)};$$

the constant suppressed here is independent of $\psi$ and $r$. Since $\psi \in C_c(\Omega_1)$ is arbitrary, this, together with (5.3), yields

$$\|g_{r^2}\|_{F_\alpha^1} \approx \|g_{r^2}\|_{\mathcal{L}_\alpha^{1,1}} \leq Cr^{-2|\beta|}\|\nu\|$$

for some constant $C > 0$ independent of $r$. So, we conclude (5.11) by Fatou's lemma. Finally, since $g \in F_\alpha^1$, one may verify $\nu = \langle \, \cdot \, , g \rangle_s$, as before. This completes the proof. $\qquad \square$

As an immediate consequence of Theorems 4.2, 5.5 and 5.7, we obtain the following duality for the Fock-Sobolev spaces.

$\boxed{\text{FSdual}}$ **Corollary 5.9.** *Let $s$, $t$, $\alpha$, $\beta$ be real numbers. Then the following dualities hold with equivalent norms:*

(1) *If $1 < p < \infty$ and $q$ is the conjugate index of $p$, then*

$$(F_{\alpha,s}^p)^* = F_{\beta,t}^q$$

*under the pairing $\langle \, , \rangle_\sigma$ where $\sigma = \alpha/p + \beta/q - s - t$;*

(2) *If $0 < p \leq 1$, then*

$$(F_{\alpha,s}^p)^* = F_\beta^\infty$$

*under the pairing $\langle \, , \rangle_\sigma$ where $\sigma = \alpha/p + \beta - s$.*

### 5.3. Complex interpolation.

We briefly recall the notion of the complex interpolation. Let $X_0$ and $X_1$ be Banach spaces both continuously imbedded in a Hausdorff topological vector space. The space $X_0 + X_1$ is a Banach space with the standard norm. Writing $S$ for the open strip consisting of complex numbers with real parts between 0 and 1, we denote by $\mathscr{F}(X_0, X_1)$ the class of all functions $\Phi : \overline{S} \to X_0 + X_1$ satisfying the following properties:

(a) $\Phi$ is holomorphic on $S$;
(b) $\Phi$ is continuous and bounded $\overline{S}$;
(c) For $j = 0, 1$ the map $x \mapsto \Phi(j + ix)$ is continuous from $\mathbb{R}$ into $X_j$.

Equipped with the norm

$$\|\Phi\|_{\mathscr{F}} := \max\left(\sup_{x \in \mathbb{R}} \|\Phi(ix)\|_{X_0}, \ \sup_{x \in \mathbb{R}} \|\Phi(1 + ix)\|_{X_1}\right),$$

the space $\mathscr{F}(X_0, X_1)$ is a Banach space. For $0 \leq \theta \leq 1$, define $[X_0, X_1]_\theta$ to be the space of all $u = \Phi(\theta)$ for some $\Phi \in \mathscr{F}(X_0, X_1)$. The norm

$$\|u\|_\theta := \inf\{\|\Phi\|_{\mathscr{F}} : u = \Phi(\theta), \Phi \in \mathscr{F}(X_0, X_1)\}$$

turns $[X_0, X_1]_\theta$ into a Banach space. As is well known, $[X_0, X_1]_\theta$ is an interpolation space between $X_0$ and $X_1$. For details we refer to [10, Chapter 2].

We also recall the well-known result concerning the complex interpolation of weighted $L^p$-spaces. Recall that, given $1 \leq p \leq \infty$ and a positive weight $\omega$ on $\mathbb{C}^n$, the $\omega$-weighted Lebesgue space $L_\omega^p = L_\omega^p(dV)$ is the space consisting of all Lebesgue measurable functions $\psi$ on $\mathbb{C}^n$ such that $\psi\omega \in L^p(dV)$. The norm of $\psi \in L_\omega^p$ is given by $\|\psi\|_{L_\omega^p} := \|\psi\omega\|_{L^p(dV)}$. For example, for the weight function defined by $\omega_r(z) = e^{\frac{1}{2}|z|^2}|z|^{-\alpha}$ for $|z| \geq r$ and $\omega_r(z) = 0$ for $|z| < r$, we have $L_{\omega_r}^p = \mathcal{L}_{\alpha p}^{p,r}$ for $1 \leq p < \infty$ but $L_{\omega_r}^\infty = \mathcal{L}_\alpha^{\infty,r}$.

In the next lemma, which is a special case of the Stein interpolation theorem, the notation $\mathcal{L}_{\alpha p}^{p,r}$ with $p = \infty$ stands for $\mathcal{L}_\alpha^{\infty,r}$.

**Lemma 5.10.** *Let $\alpha_0$, $\alpha_1$ be real numbers, $1 \leq p_0 \leq p_1 \leq \infty$ and $r > 0$. Let $0 \leq \theta \leq 1$. Then*

$$[\mathcal{L}_{\alpha_0 p_0}^{p_0,r}, \mathcal{L}_{\alpha_1 p_1}^{p_1,r}]_\theta = \mathcal{L}_{\alpha p}^{p,r}$$

*with equal norms where*

$$\frac{1}{p} = \frac{1-\theta}{p_0} + \frac{\theta}{p_1} \ \ and \ \ \alpha = (1-\theta)\alpha_0 + \theta\alpha_1. \tag{5.14}$$

In the next theorem we use the notation $F_{\alpha p}^p$ with $p = \infty$ for the space $F_\alpha^\infty$ for a unified statement.

**Theorem 5.11.** *Let $\alpha_0$, $\alpha_1$ be real numbers and $1 \leq p_0 \leq p_1 \leq \infty$. Let $0 \leq \theta \leq 1$. Then*

$$[F_{\alpha_0 p_0}^{p_0}, F_{\alpha_1 p_1}^{p_1}]_\theta = F_{\alpha p}^p$$

*with equivalent norms where $p$ and $\alpha$ are as in* (5.14).

*Proof.* Using Theorem 5.4, fix an $r > 0$ sufficiently small so that $P_\alpha^r : F_\alpha^p \to F_\alpha^p$ is invertible. We see from (5.3), Lemma 5.10 and the definition of complex interpolation that $[F_{\alpha_0 p_0}^{p_0}, F_{\alpha_1 p_1}^{p_1}]_\theta$ is continuously embedded into $F_\alpha^p$.

Conversely, let $f \in F_\alpha^p$ and pick $g \in F_\alpha^p$ such that $P_\alpha^r g = f$. By Lemma 5.10 we can pick some function $\Phi \in \mathscr{F}(\mathcal{L}_{\alpha_0 p_0}^{p_0, r}, \mathcal{L}_{\alpha_1 p_1}^{p_1, r})$ and a constant $C > 0$ such that $\Phi(\theta) = g$ and

$$\sup_{\operatorname{Re}\lambda=j} \|\Phi(\lambda)\|_{\mathcal{L}_{\alpha_j p_j}^{p_j, r}} \leq C\|g\|_{F_{\alpha p}^p} \tag{5.15}$$ `phi`

for each $j = 0, 1$. Since $P_\alpha^r : \mathcal{L}_{\alpha_j p_j}^{r} \to F_{\alpha_j p_j}^{p_j}$ is bounded for each $j$ by Lemma 5.2, we have by (5.15)

$$\sup_{\operatorname{Re}\lambda=j} \left\| P_\alpha^r[\Phi(\lambda)] \right\|_{F_{\alpha_j p_j}^{p_j}} \lesssim \|g\|_{F_{\alpha p}^p}$$

for each $j$. Now, defining

$$\Psi(\lambda) := P_\alpha^r[\Phi(\lambda)], \qquad 0 \leq \operatorname{Re}\lambda \leq 1,$$

we see that $\Psi \in \mathscr{F}(F_{\alpha_0 p_0}^{p_0}, F_{\alpha_1 p_1}^{p_1})$. Moreover, we have $\Psi(\theta) = P_\alpha^r[\Phi(\theta)] = P_\alpha^r g = f$. So, we have $f \in [F_{\alpha_0 p_0}^{p_0}, F_{\alpha_1 p_1}^{p_1}]_\theta$ with $\|f\|_\theta \lesssim \|f\|_{F_{\alpha p}^p}$. Thus we conclude that $F_\alpha^p$ is continuously embedded into $[F_{\alpha_0 p_0}^{p_0}, F_{\alpha_1 p_1}^{p_1}]_\theta$. The proof is complete. $\square$

As a consequence of Theorems 4.2 and 5.11, we obtain the following complex interpolation for the Fock-Sobolev spaces. In the next corollary we also use the notation $F_{\alpha p, s}^p$ with $p = \infty$ for the space $F_{\alpha, s}^\infty$ for a unified statement.

`FSWcor` **Corollary 5.12.** *Let* $\alpha_j$, $s_j$ *be real numbers for* $j = 0, 1$ *and* $1 \leq p_0 \leq p_1 \leq \infty$. *Let* $0 \leq \theta \leq 1$. *Then*

$$[F_{\alpha_0 p_0, s_0}^{p_0}, F_{\alpha_1 p_1, s_1}^{p_1}]_\theta = F_{\alpha p, s}^p$$

*with equivalent norms where* $p$, $\alpha$ *are as in* (5.14) *and* $s = (1-\theta)s_0 + \theta s_1$.

5.4. **Carleson measure.** Let $0 < p < \infty$ and $\alpha$ be real. Let $\mu$ be a positive Borel measure on $\mathbb{C}^n$. We say that $\mu$ is a *Carleson measure for* $F_\alpha^p$ if there is a constant $C > 0$ such that

$$\int_{\mathbb{C}^n} \left| f(z) e^{-\frac{1}{2}|z|^2} \right|^p d\mu(z) \leq C\|f\|_{F_\alpha^p}^p \tag{5.16}$$ `t:carleson`

for all $f \in F_\alpha^p$. We say that $\mu$ is a *vanishing Carleson measure for* $F_\alpha^p$ if

$$\lim_{j \to \infty} \int_{\mathbb{C}^n} \left| f_j(z) e^{-\frac{1}{2}|z|^2} \right|^p d\mu(z) = 0$$

for every bounded sequence $\{f_j\}$ in $F_\alpha^p$ that converges to $0$ uniformly on compact subsets of $\mathbb{C}^n$.

In what follows $B(z, r)$ denotes the Euclidean ball centered at $z \in \mathbb{C}^n$ with radius $r > 0$.

`t:Carleson` **Theorem 5.13.** *Let* $0 < p < \infty$, $r > 0$, *and* $\alpha$ *be a real number. Let* $\mu$ *be a positive Borel measure on* $\mathbb{C}^n$. *Then the following statements hold:*

(1) $\mu$ is a Carleson measure for $F_\alpha^p$ if and only if there is a constant $C > 0$ such that

$$\mu[B(z,r)] \leq \frac{C}{(1+|z|)^\alpha} \qquad (5.17)$$ `t:measure`

for all $z \in \mathbb{C}^n$.

(2) $\mu$ is a vanishing Carleson measure for $F_\alpha^p$ if and only if

$$\mu[B(z,r)](1+|z|)^\alpha \to 0 \qquad (5.18)$$ `t:vanish`

as $|z| \to \infty$.

*Proof.* Using the Fock kernels as test functions and utilizing Lemma 3.6, one may imitate the proof of [5, Theorem 21] to prove the necessities of (1) and (2).

For the sufficiencies, while one may also imitate the proof of [5, Theorem 21], we provide simpler proofs. First, assume (5.17) and fix $f \in F_\alpha^p$. By Lemma 2.1, there is a constant $C = C(p,r) > 0$ such that

$$|f(z)|^p e^{-\frac{p}{2}|z|^2} \leq C \int_{B(z,r/2)} \left| f(w) e^{-\frac{1}{2}|w|^2} \right|^p \, dV(w). \qquad (5.19)$$ `subine`

for all $z \in \mathbb{C}^n$. Note $z \in B(w,r)$ for every $w \in B(z,r/2)$. Thus, integrating both sides of the above against the measure $\mu$ and then interchanging the order of integrations, we have

$$\int_{\mathbb{C}^n} \left| f(z) e^{-\frac{1}{2}|z|^2} \right|^p \, d\mu(z) \lesssim \int_{\mathbb{C}^n} \left\{ \int_{B(z,r/2)} \left| f(w) e^{-\frac{1}{2}|w|^2} \right|^p \, dV(w) \right\} d\mu(z)$$
$$\leq \int_{\mathbb{C}^n} \left| f(w) e^{-\frac{1}{2}|w|^2} \right|^p \mu[B(w,r)] \, dV(w),$$

which, together with (5.17), yields (5.16). This completes the proof of the sufficiency of (1).

Next, assume (5.18). Let $\{f_j\}$ be a bounded sequence in $F_\alpha^p$ that converges to $0$ uniformly on compact subsets of $\mathbb{C}^n$ and put

$$I(f_j) := \int_{\mathbb{C}^n} \left| f_j(z) e^{-\frac{1}{2}|z|^2} \right|^p \, d\mu(z)$$

for each $j$. Fix an arbitrary $R > 0$. Since $\{f_j\}$ converges uniformly on compact subsets of $\mathbb{C}^n$, we have

$$\limsup_{j \to \infty} I(f_j) \leq \limsup_{j \to \infty} \int_{|z| \geq R+r/2} \left| f_j(z) e^{-\frac{1}{2}|z|^2} \right|^p \, d\mu(z).$$

Meanwhile, proceeding as above, we have by (5.19)

$$\int_{|z|\geq R+r/2} \left| f_j(z) e^{-\frac{1}{2}|z|^2} \right|^p \, d\mu(z)$$

$$\lesssim \int_{|z|\geq R+r/2} \left\{ \int_{B(z,r/2)} \left| f_j(w) e^{-\frac{1}{2}|w|^2} \right|^p \, dV(w) \right\} d\mu(z)$$

$$\lesssim \int_{|w|\geq R} \mu[B(w,r)] \left| f_j(w) e^{-\frac{1}{2}|w|^2} \right|^p \, dV(w)$$

$$\lesssim M \sup_{|w|\geq R} \mu[B(w,r)](1+|w|)^\alpha$$

where $M := \sup_j \|f_j\|_{F_\alpha^p}^p < \infty$. So, we obtain

$$\limsup_{j\to\infty} I(f_j) \leq CM \sup_{|w|\geq R} \mu[B(w,r)](1+|w|)^\alpha$$

for some constant $C > 0$ independent of $R > 0$. Now, taking $R \to \infty$, we conclude by (5.18) that $I(f_j) \to 0$ as $j \to \infty$. This completes the proof of the sufficiency of (2). The proof is complete. $\qquad\square$

DEPARTMENT OF MATHEMATICS, PUSAN NATIONAL UNIVERSITY, PUSAN 609-735, REPUBLIC OF KOREA
*E-mail address*: chohr@pusan.ac.kr

DEPARTMENT OF MATHEMATICS, KOREA UNIVERSITY, SEOUL 136-713, REPUBLIC OF KOREA
*E-mail address*: cbr@korea.ac.kr

DEPARTMENT OF MATHEMATICS, KOREA UNIVERSITY, SEOUL 136-713, REPUBLIC OF KOREA
*E-mail address*: koohw@korea.ac.kr